\documentclass[12pt]{imsart}
\usepackage{amsthm,amsmath}
\usepackage[colorlinks,citecolor=blue,urlcolor=blue]{hyperref}

\usepackage{amsmath,amsfonts,amssymb,amsthm}

\usepackage{enumitem}

\RequirePackage{amsthm,amsmath,amsfonts,amssymb}
\RequirePackage[numbers]{natbib}
\RequirePackage[colorlinks,citecolor=blue,urlcolor=blue]{hyperref}
\RequirePackage{graphicx}
\usepackage{lineno,hyperref}

\startlocaldefs
\modulolinenumbers[5]

\usepackage[utf8]{inputenc}
\usepackage[english]{babel}

\usepackage[margin=1in]{geometry}
\usepackage{bm}
\usepackage{graphicx}
\usepackage{amssymb,amsmath,amsthm,mathrsfs}
\usepackage{epstopdf}
\usepackage{algorithm}
\usepackage{amsfonts,delarray}
\usepackage{algorithm,algcompatible,amsmath}
\usepackage{rotating,color}
\usepackage{subfigure}
\usepackage{multirow}
\usepackage{titlesec,textcase}
\usepackage{longtable}
\usepackage{bbm}
\usepackage{rotating}

\newtheorem{Theorem}{Theorem}[section]
\newtheorem{Lemma}[Theorem]{Lemma}
\newtheorem{Corollary}[Theorem]{Corollary}
\newtheorem{Proposition}[Theorem]{Proposition}

\newtheorem{Definition}[Theorem]{Definition}

\theoremstyle{definition}

\RequirePackage[normalem]{ulem} 

\definecolor{rp}{RGB}{83,54,106}

\def\boxit#1{\vbox{\hrule\hbox{\vrule\kern6pt\vbox{\kern6pt#1\kern6pt}\kern6pt\vrule}\hrule}}

\endlocaldefs

\allowdisplaybreaks

\begin{document}
\begin{frontmatter}
\title{Asymptotic distributions of the average clustering coefficient and its variant}

\runtitle{Asymptotic distribution of the clustering coefficient }
\runauthor{ }
\begin{aug}

\author[A]{\fnms{Mingao} \snm{Yuan}\ead[label=e1]{mingao.yuan@ndsu.edu}},
\author[C]{\fnms{Xiaofeng} \snm{Zhao}\ead[label=e3]{zxfstats@ncwu.edu.cn}}



\address[A]{Department of Statistics,
North Dakota State University,
\printead{e1}}

\address[C]{School of Mathematics and Statistics,
North China University of Water Resources and Electric Power,
\printead{e3}}
\end{aug}

\begin{abstract}
In network data analysis, summary statistics of a network can provide us with meaningful insight into the structure of the network. The average clustering coefficient is one of the most popular and widely used network statistics.
In this paper, we investigate the asymptotic distributions of the average clustering coefficient and its variant of a heterogeneous  Erd\"{o}s-R\'{e}nyi random graph. We show that the standardized average clustering coefficient converges in distribution to the standard normal distribution. Interestingly, the variance of the average clustering coefficient exhibits a
phase transition phenomenon. The sum of weighted triangles is a variant of the average clustering coefficient. It is recently introduced to detect geometry in a network. We also derive the asymptotic distribution of the sum weighted triangles, which does not exhibit a
phase transition phenomenon as the average clustering coefficient. This result signifies the difference between the two summary statistics.
\end{abstract}



\end{frontmatter}

\section{Introduction}
\label{S:1}

Network data consists of a set of vertices and a set of edges that represents the connections between vertices. Network data analysis is widely used to study many real-world problems. Amazon and eBay utilize  the recommendation network to advertise their online goods \cite{LSK1}. The telecommunication companies optimize the performance of 5G wireless network by studying its topology structure \cite{ARS01}. The academic journals investigate the citation network to study the relationship between papers, authors, and scientific work \cite{RFV1}. Biological network is used to detect gene-gene interactions \cite{CV1}. Due to the widespread applications, network analysis becomes one of the most proactive research directions.   

In network data analysis, a typical task is to
understand the structural properties of a network. Visualization is perhaps the most  straightforward method to describe the structure of a network. However, larger networks can be difficult to envision and describe. Numerous summary statistics have therefore been proposed to quantify the structures of a network. Based on these statistics, we are able to compare networks or classify
them according to properties that they exhibit.

One commonly used network statistic is the average clustering coefficient, which measures the trend of the vertices of a network to cluster together \cite{WS98}. The average clustering coefficient was first introduced in social network analysis to quantify the property that friends of a
friend tend to be friends \cite{BW00}.   Since then, the average clustering coefficient has been widely used in network
analysis. \cite{XKB1} used the average clustering coefficient to  identify functional disconnections in patients with borderline personality disorder. \cite{CLX01} applied the average clustering coefficient to examine the dissimilarities of the brain functional networks between endurance runners and healthy controls. \cite{TPL01} applied the average clustering coefficient to investigate the functional brain network development in children. 
In \cite{BBA01}, the average clustering coefficient was employed to compare the brain connectivity networks with
generative models of structured networks. In \cite{TOSHK06}, the clustering coefficient was used to quantify the network structure of large networks.

It is shown in \cite{MLS22} that the average clustering coefficient is insufficient to signal the presence of hyperbolic geometry in a network. Therefore, \cite{MLS22} introduces the sum of weighted triangles to detect the geometry in a network. Note that the 
 average clustering coefficient is also a sum of weighted triangles. However, the weights are different from the sum of weighted triangles introduced in \cite{MLS22}. Hence, the sum of weighted triangles is a variant of the average clustering coefficient. The analytical analysis and numeric studies in  \cite{MLS22} show that the sum of weighted triangles  has high potential for uncovering a hidden
network geometry.

One of the important research topics in network analysis is to study the asymptotic properties of summary statistics \cite{BM06,CHHS20,Y23,Y23b,Y23c, ZY23}. In this paper, we are interested in the asymptotic distributions of the average clustering coefficient and the sum of weighted triangles in a heterogeneous  Erd\"{o}s-R\'{e}nyi random graph. Note that both the average clustering coefficient and the sum of weighted triangles are sums of dependent terms. It is not a trivial task to derive their asymptotic distributions. We prove that the standardized average clustering coefficient and the standardized sum of weighted triangles converge in distribution to the standard normal distribution. Interestingly, the variance of the average clustering coefficient exhibits a phase transition phenomenon. However, the sum of weighted triangles does not present similar phase transition. These results highlight the difference between the two summary statistics.

The rest of the paper is organized as follows. In section 2, we introduce the heterogeneous  Erd\"{o}s-R\'{e}nyi random graph model, definitions of the average clustering coefficient and the sum of the weighted triangles, and present the main results. The proof is deferred to section 3.

\medskip

\noindent
{\bf Notation:} We adopt the  Bachmann–Landau notation throughout this paper. Let $a_n$  and $b_n$ be two positive sequences. Denote $a_n=\Theta(b_n)$ if $c_1b_n\leq a_n\leq c_2 b_n$ for some positive constants $c_1,c_2$. Denote  $a_n=\omega(b_n)$ if $\lim_{n\rightarrow\infty}\frac{a_n}{b_n}=\infty$. Denote $a_n=O(b_n)$ if $a_n\leq cb_n$ for some positive constants $c$. Denote $a_n=o(b_n)$ if $\lim_{n\rightarrow\infty}\frac{a_n}{b_n}=0$. Let $\mathcal{N}(0,1)$ be the standard normal distribution and $X_n$ be a sequence of random variables. Then $X_n\Rightarrow\mathcal{N}(0,1)$ means $X_n$ converges in distribution to the standard normal distribution as $n$ goes to infinity. Denote $X_n=O_P(a_n)$ if $\frac{X_n}{a_n}$ is bounded in probability. Denote $X_n=o_P(a_n)$ if $\frac{X_n}{a_n}$ converges to zero in probability as $n$ goes to infinity. Let $\mathbb{E}[X_n]$ and $Var(X_n)$ denote the expectation and variance of a random variable $X_n$ respectively. $\mathbb{P}[E]$ denote the probability of an event $E$. Let $f=f(x)$ be a function. Denote $f^{(k)}(x)=\frac{d^kf}{dx^k}(x)$ for any positive integer $k$. $\exp(x)$ denote the exponential function $e^x$.
For positive integer $n$, denote $[n]=\{1,2,\dots,n\}$. Given a finite set $E$, $|E|$ represents the number of elements in $E$. Given positive integer $t$, $\sum_{i_1\neq i_2\neq\dots\neq i_t}$ means summation over all integers $i_1,i_2,\dots,i_t$ in $[n]$ such that $|\{i_1,i_2,\dots,i_t\}|=t$. $\sum_{i_1< i_2<\dots< i_t}$ means summation over all integers $i_1,i_2,\dots,i_t$ in $[n]$ such that $i_1<i_2<\dots<i_t$.

\section{Main results}

A graph is a mathematical model that consists of a set of nodes (vertices) and a set of edges. Let $\mathcal{V}=[n]$ for any positive integer $n$.
An \textit{undirected} graph on $\mathcal{V}$ is defined as the pair $\mathcal{G}=(\mathcal{V},\mathcal{E})$, where $\mathcal{E}$ is a set of subsets of $\mathcal{V}$ such that $|e|=2$ for every $e\in\mathcal{E}$. Each element in $\mathcal{V}$ is called a node or vertex of the graph and each element in $\mathcal{E}$ is called an edge. A graph can be conveniently represented as an adjacency matrix $A$. In $A$, $A_{ij}=1$ if $\{i,j\}$ is an edge, $A_{ij}=0$ otherwise and $A_{ii}=0$. Since $\mathcal{G}$ is undirected, the adjacency matrix $A$ is symmetric. The degree $d_i$ of node $i$ is the number of edges connecting it, that is, $d_i=\sum_{j}A_{ij}$.  A graph is said to be random if $A_{ij} (1\leq i<j\leq n)$ are random.

\begin{Definition}\label{def1}
Let $\alpha$ and $\beta$ be constants between zero and one, that is, $\alpha,\beta\in(0,1)$, and 
\[W=\{w_{ij}\in[\beta,1]| 1\leq i,j\leq n, w_{ji}=w_{ij} ,w_{ii}=0\}.\]
Define a heterogeneous random graph  $\mathcal{G}_n(\alpha,\beta, W)$ as  
\[\mathbb{P}(A_{ij}=1)=p_nw_{ij},\]
where $A_{ij}$ $(1\leq i<j\leq n)$ are independent, $A_{ij}=A_{ji}$ and $p_n=n^{-\alpha}$.
\end{Definition}

The random graph  $\mathcal{G}_n(\alpha,\beta, W)$ is heterogeneous because the expected degrees of nodes may be different. Specifically, $\mathbb{E}[d_i]=\sum_{k}p_nw_{ik}$. In general, $\mathbb{E}[d_i]\neq \mathbb{E}[d_j]$ for distinct nodes $i,j$.
If $w_{ij}=c$ $(1\leq i<j\leq n)$ for a constant $c\in(0,1)$, the expected degrees of nodes are the same. In this case, $\mathcal{G}_n(\alpha,\beta, W)$ is homogeneous and it is called the Erd\"{o}s-R\'{e}nyi random graph. For convenience, we simply denote it as $\mathcal{G}_n(\alpha)$.  The random graph  $\mathcal{G}_n(\alpha,\beta, W)$ serves as our benchmark model. It is studied in \cite{Y23,Y23c} and includes the inhomogeneous Erd\"{o}s-R\'{e}nyi random graphs in  \cite{CHHS20,Y23b,ZY23} as a special case.

\subsection{The average clustering coefficient}
In social networks, nodes tend to develop highly connected
neighborhoods \cite{WS98,BW00}. The average
clustering coefficient was first introduced to measure such properties of social networks \cite{WS98,BW00}.
The clustering coefficient of each node in a graph is the fraction of triangles that actually exist over all possible triangles in its neighborhood. It measures the triangular pattern and the connectivity in a node’s neighborhood. A node has a high clustering
coefficient if its neighbors tend to be directly connected with each other.  The average clustering coefficient of a graph is the mean of the clustering coefficients of all nodes. It measures the extent to which nodes in a graph tend to cluster together.  In terms of the adjacency matrix $A$, the average clustering coefficient is defined as
\begin{equation}\label{clustercoe}
\overline{\mathcal{C}}_n=\frac{1}{n}\sum_{i=1}^n\frac{\sum\limits_{j\neq k} A_{ij}A_{jk}A_{ki}}{\sum\limits_{j\neq k}A_{ij}A_{ik}},
\end{equation}
where any summation terms with denominator zero is set to be zero. 

Recently, \cite{NS21} studied the robustness of the average clustering coefficient and \cite{ZY23} obtained the limit of the average clustering coefficient. In this paper, we derive its asymptotic distribution as follows.

\begin{Theorem}\label{mainthm}
For the heterogeneous random graph $\mathcal{G}_n(\alpha,\beta, W)$, we have
    \[\frac{\overline{\mathcal{C}}_n-\mathbb{E}[\overline{\mathcal{C}}_n]}{\sigma_n}\Rightarrow \mathcal{N}(0,1),\]
where $\sigma_n^2=\sigma_{1n}^2+\sigma_{2n}^2$, $\mu_i=\sum_{j}\mu_{ij}$, $t_i=\sum\limits_{ j\neq k} A_{ij}A_{jk}A_{ki}$,  and
\[\sigma_{1n}^2=\frac{4}{n^2}\sum_{i<j<k}a_{ijk}^2\mu_{ij}(1-\mu_{ij})\mu_{jk}(1-\mu_{jk})\mu_{ki}(1-\mu_{ki}),\ \ \ \ \sigma_{2n}^2=\frac{4}{n^2}\sum_{i<j}e_{ij}^2\mu_{ij}(1-\mu_{ij}).\]
\[a_i=\mathbb{E}\left[\frac{1}{d_i(d_i-1)}\right],\ \ \ b_i=\frac{\mathbb{E}[t_i](2\mu_i-1)}{\mu_i^2(\mu_i-1)^2},\ \ c_{ij}=\sum_{k\neq i,j}a_k\mu_{ki}\mu_{kj},\ \ \]
\[e_{ij}=2c_{ij}+2(a_id_{ij}+a_jd_{ji})-b_i-b_j,\ \ \ \ a_{ijk}=a_i+a_j+a_k,\ \ \ d_{ij}=\sum_{k\neq i,j}\mu_{ki}\mu_{kj}.\]
\end{Theorem}

Based on Theorem \ref{mainthm}, the standardized  average clustering coefficient  converges in distribution to the standard normal distribution. 
The proof of Theorem \ref{mainthm} is not trivial, due to the fact that $n\overline{\mathcal{C}}_n$ is a sum of dependent random variables. Our proof strategy is as follows: use Taylor expansion to  expand the summation terms to $k_0=\lceil 1+\frac{1}{1-\alpha}\rceil +2$ order, isolate the leading term and prove the leading term converges in distribution to the standard normal distribution.

Next we apply  Theorem \ref{mainthm} to 
the Erd\"{o}s-R\'{e}nyi random graph $\mathcal{G}_n(\alpha)$ and get the following corollary.

\begin{Corollary}\label{mcor}
\label{mainthmcor}
     For the Erd\"{o}s-R\'{e}nyi random graph $\mathcal{G}_n(\alpha)$, we have
    \[\frac{\overline{\mathcal{C}}_n-\mathbb{E}[\overline{\mathcal{C}}_n]}{\sigma_n}\Rightarrow \mathcal{N}(0,1),\]
    where $\sigma_n^2=\sigma_{1n}^2+\sigma_{2n}^2$, 
    $\sigma_{1n}^2=\frac{6}{n^{3-\alpha}}(1+o(1))$ and $\sigma_{2n}^2=\frac{2}{n^{2+\alpha}}(1+o(1))$. Hence we have
\begin{equation*} 
\sigma_n^2 =
    \begin{cases}
         \frac{6}{n^{3-\alpha}}(1+o(1)), & \text{if $\alpha>\frac{1}{2}$}, \\
   \frac{8}{n^2\sqrt{n}}(1+o(1)), & \text{if $\alpha=\frac{1}{2}$}, \\
    \frac{2}{n^{2+\alpha}}(1+o(1)) ,& \text{if $\alpha<\frac{1}{2}$}.     
    \end{cases}       
\end{equation*}
\end{Corollary}
According to Corollary \ref{mcor}, the average clustering coefficient of the Erd\"{o}s-R\'{e}nyi random graph $\mathcal{G}_n(\alpha)$ shows a phase change phenomenon. Given large $n$, we have
\[\lim_{\alpha \to (\frac{1}{2})^{-}}\frac{2}{n^{2+\alpha}}=\frac{2}{n^{2}\sqrt{n}}\neq \frac{8}{n^2\sqrt{n}}, \]

\[\lim_{\alpha \to (\frac{1}{2})^{+}}\frac{6}{n^{3-\alpha}}=\frac{6}{n^{2}\sqrt{n}}\neq \frac{8}{n^2\sqrt{n}}. \]
As $\alpha$ varies around $\frac{1}{2}$, $\sigma_n^2$ does not change continuously as a function of $\alpha$. In this sense, the scale $\sigma_n^2$ exhibits a phase change phenomenon at $\alpha=\frac{1}{2}$.

\subsection{The sum of weighted triangles}

Many real-world networks have geometric structures and can be accurately modelled by geometric random graphs. Nevertheless, the presence of geometry in an observed network is not always evident. 
The average clustering coefficient is a standard statistic to indicate the presence of geometry. However, it is shown that the average clustering coefficient fails to detect hyperbolic geometry in a network \cite{MLS22}. Therefore, \cite{MLS22} introduces a novel triangle-based
statistic, that is, the sum of weighted triangles. Note that the
 average clustering coefficient can also be written as a sum of weighted triangles, with the weights different from the sum of weighted triangles introduced in \cite{MLS22}. Hence, the sum of weighted triangles is a variant of the average clustering coefficient. By the analytical analysis and numeric studies in  \cite{MLS22}, the sum of weighted triangles shows high potential for uncovering hidden
network geometry.  

In terms of the adjacency matrix $A$, the sum of the weighted triangles is defined as \cite{MLS22}
\begin{eqnarray*}
\mathcal{T}_n=\sum_{i<j<k}\frac{A_{ij}A_{jk}A_{ki}}{d_id_jd_k},
\end{eqnarray*}
where the summation term is set to be zero if $d_id_jd_k=0$.  \cite{MLS22} obtains bounds of $\mathcal{T}_n$ in the inhomogeneous random graph and the geometric inhomogeneous random graph. In this paper, we derive the asymptotic distribution of  $\mathcal{T}_n$ in the heterogeneous random graph $\mathcal{G}_n(\alpha,\beta, W)$.

\begin{Theorem}\label{mainthm2}
    For the heterogeneous random graph $\mathcal{G}_n(\alpha,\beta, W)$, we have
    \[\frac{\mathcal{T}_n-\mathbb{E}[\mathcal{T}_n]}{v_n}\Rightarrow \mathcal{N}(0,1),\]
 where $v_n^2=v_{1n}^2+v_{2n}^2$ and
\[v_{1n}^2=\sum_{i<j<k}\frac{\mu_{ij}(1-\mu_{ij})\mu_{jk}(1-\mu_{jk})\mu_{ki}(1-\mu_{ki})}{\mu_i^2\mu_j^2\mu_k^2},\]
\[v_{2n}^2=\sum_{i<j}\left(\gamma_{ij}-\frac{\eta_i+\eta_j}{2}\right)^2\mu_{ij}(1-\mu_{ij}),\]
 \[\eta_{i}=\sum_{\substack{j\neq k}}\frac{\mu_{ij}\mu_{jk}\mu_{ki}}{\mu_i^2\mu_j\mu_k}, \ \ \ \ \  \ \ \gamma_{ij}=\sum_{k\not\in\{i,j\}}\frac{\mu_{jk}\mu_{ki}}{\mu_i\mu_j\mu_k}.\]
\end{Theorem}

Based on Theorem \ref{mainthm2}, the standardized sum of weighted triangles converges in distribution to the standard normal distribution.

 Next we apply Theorem \ref{mainthm2} to a special heterogeneous random graph.
In the definition of $\mathcal{G}_n(\alpha,\beta,W)$, let 
$w_{ij}=w_iw_j$ for $w_i\in[\beta,1]$. In this case, $\mathcal{G}_n(\alpha,\beta,W)$ is called rank-1 heterogeneous random graph.

\begin{Corollary}\label{rank1cor}
    For the rank-1 heterogeneous random graph  $\mathcal{G}_n(\alpha,\beta,W)$, we have
    \[\frac{\mathcal{T}_n-\mathbb{E}[\mathcal{T}_n]}{v_n}\Rightarrow \mathcal{N}(0,1),\]
    where $v_{n}^2=\frac{n^3}{6w^6p_n^3}(1+o(1))$ and $w=\sum_{i}w_i$. Especially, for the Erd\"{o}s-R\'{e}nyi random graph $\mathcal{G}_n(\alpha)$,
    $v_{n}^2=\frac{1}{6n^{3(1-\alpha)}}(1+o(1))$.
\end{Corollary}

Note that
\[\frac{n^3}{6w^6p_n^3}=\frac{1}{6\left(\frac{w}{n}\right)^6n^{3(1-\alpha)}},\]
and $\frac{w}{n}$ is independent of $\alpha$. As $\alpha$ varies, $v_n^2$ changes continuously as a function of $\alpha$. Hence, the sum of weighted triangles of the rank-1 random graph  $\mathcal{G}_n(\alpha,\beta,W)$ does not exhibit a phase change phenomenon. This signifies the difference between the sum of weighted triangles and the average clustering coefficient.

\section{Proof of main results}

In this section, we provide the detailed proofs of the main results. 
Denote $\bar{A}_{ij}=A_{ij}-\mu_{ij}$ in this section.
Before proving Theorem \ref{mainthm} and Theorem \ref{mainthm2},
we present two lemmas.  

\begin{Lemma}\label{lem1}\
Let $\delta_n=\left(\log(np_n)\right)^{-2}$. 
For the random graph $\mathcal{G}_n(\alpha,\beta, W)$,  
 we have
\[\mathbb{P}(d_{i}=k)\leq \exp(-np_n\beta(1+o(1))),\ \ \ \ k\leq \delta_nnp_n,\]
uniformly for all $i\in[n]$.
\end{Lemma}
Lemma \ref{lem1} is proved in \cite{Y23b}. Hence we omit the proof.

\begin{Lemma}\label{lem2}
 For the random graph $\mathcal{G}_n(\alpha,\beta, W)$ and an even positive integer $s$,  we have
\[\mathbb{E}[(d_i-\mu_i)^s]=O\left((np_n)^{\frac{s}{2}}\right).\]
\end{Lemma}

\noindent
{\bf Proof of Lemma \ref{lem2}.} Given integer $t\leq s$, let $\lambda_{t1},\lambda_{t2},\dots,\lambda_{tt}$ be arbitrary non-negative integers such that $\lambda_{t1}+\lambda_{t2}+\dots+\lambda_{tt}=s$. Then
\begin{eqnarray}\nonumber
(d_i-\mu_i)^s&=&\sum_{j_1,\dots,j_s\neq i}\bar{A}_{ij_1}\bar{A}_{ij_2}\dots \bar{A}_{ij_s}\\ \label{dseq}
&=&\sum_{t=1}^s\sum_{\lambda_{t1},\dots,\lambda_{tt}}\sum_{j_1\neq  \dots\neq j_t\neq i}\bar{A}_{ij_1}^{\lambda_{t1}}\bar{A}_{ij_2}^{\lambda_{t2}}\dots \bar{A}_{ij_t}^{\lambda_{tt}}.
\end{eqnarray}
It is easy to verify that $\mathbb{E}[\bar{A}_{ij}^k]=0$ if $k=1$ and $\mathbb{E}[\bar{A}_{ij}^k]=O(p_n)$ if $k\geq2$. For distinct indices $j_1,j_2,\dots,j_t$, $\bar{A}_{ij_1},\bar{A}_{ij_2},\dots,\bar{A}_{ij_t}$  are independent. Hence, we have
\begin{equation*} 
\mathbb{E}[\bar{A}_{ij_1}^{\lambda_{t1}}\bar{A}_{ij_2}^{\lambda_{t2}}\dots \bar{A}_{ij_t}^{\lambda_{tt}}] =
    \begin{cases}
      O(p_n^{t})& \text{if $\lambda_{tl}\geq2$ for all $l$}, \\
      0 & \text{if $\lambda_{tl}=1$ for some $l$}.     
    \end{cases}       
\end{equation*}
When $\lambda_{tl}\geq2$ for all $l$, $t\leq\frac{s}{2}$. There are at most $n^t$ choices for indices $j_1,j_2,\dots,j_t$. Hence 
\begin{eqnarray}\nonumber
\sum_{j_1\neq  \dots\neq j_t}\mathbb{E}[\bar{A}_{ij_1}^{\lambda_{t1}}\bar{A}_{ij_2}^{\lambda_{t2}}\dots \bar{A}_{ij_t}^{\lambda_{tt}}]=O\left((np_n)^{\frac{s}{2}}\right).
\end{eqnarray}
Then
\[\mathbb{E}[(d_i-\mu_i)^s]=O\left((np_n)^{\frac{s}{2}}\right).\]

\qed

\begin{Lemma}\label{lem3} For the random graph $\mathcal{G}_n(\alpha,\beta, W)$, we have
\begin{eqnarray}\nonumber
\mathbb{E}[\overline{\mathcal{C}}_n]&=&\frac{1}{n}\sum_{i=1}^n\mathbb{E}[t_i]\mathbb{E}\left[\frac{1}{d_i(d_i-1)}\right]-\frac{2}{n}\sum_{i\neq j\neq k}\frac{2\mu_i-1}{\mu_i^2(\mu_i-1)^2}\mu_{ij}(1-\mu_{ij})\mu_{ik}\mu_{jk}\\
&&+o\left(\frac{1}{n\sqrt{np_n}}+\frac{\sqrt{p_n}}{n}\right).
\end{eqnarray}

\end{Lemma}

\noindent
{\bf Proof of Lemma \ref{lem3}.} Let $t_i=\sum\limits_{j\neq k}A_{ij}A_{jk}A_{ki}$ and $h(x)=\frac{1}{x(x-1)}$. The $k$-th derivative of $h(x)$ is equal to
\begin{eqnarray}\label{hdriv}
h^{(k)}(x)&=&\left(\frac{1}{x-1}-\frac{1}{x}\right)^{(k)}=k!(-1)^k\frac{\sum_{t=1}^{k+1}\binom{k+1}{t}x^{k+1-t}(-1)^{t+1}}{(x-1)^{k+1}x^{k+1}}.
\end{eqnarray}
Let $k_0=\lceil 1+\frac{1}{1-\alpha}\rceil +2$.
By Taylor expansion, we have
\begin{eqnarray}\nonumber
\frac{1}{d_i(d_i-1)}&=&\frac{1}{\mu_i(\mu_i-1)}-\frac{2\mu_i-1}{\mu_i^2(\mu_i-1)^2}(d_i-\mu_i)+\sum_{s=2}^{k_0-1}\frac{h^{(s)}(\mu_i)}{s!}(d_i-\mu_i)^s\\ \label{0eq2}
&&+\frac{h^{(k_0)}(X_i)}{k_0!}(d_i-\mu_i)^{k_0},
\end{eqnarray}
where $X_i$ is between $d_i$ and $\mu_i$. Then
\begin{eqnarray}\nonumber
\mathbb{E}[\overline{\mathcal{C}}_n]&=&\frac{1}{n}\sum_{i}\frac{\mathbb{E}[t_i]}{\mu_i(\mu_i-1)}-\frac{1}{n}\sum_{i}\frac{2\mu_i-1}{\mu_i^2(\mu_i-1)^2}\mathbb{E}[t_i(d_i-\mu_i)]\\ \label{00eq2}
&&
+\sum_{s=2}^{k_0-1}\frac{1}{n}\sum_{i}\frac{h^{(s)}(\mu_i)}{s!}\mathbb{E}[t_i(d_i-\mu_i)^s]+\frac{1}{n}\sum_{i}\mathbb{E}\left[t_i\frac{h^{(k_0)}(X_i)}{k_0!}(d_i-\mu_i)^{k_0}\right].
\end{eqnarray}

Taking expectation on both sides of (\ref{0eq2}) yields
\begin{eqnarray}\label{edd}
\mathbb{E}\left[\frac{1}{d_i(d_i-1)}\right]&=&\frac{1}{\mu_i(\mu_i-1)}+\sum_{s=2}^{k_0-1}\frac{h^{(s)}(\mu_i)}{s!}\mathbb{E}[(d_i-\mu_i)^s]+\mathbb{E}\left[\frac{h^{(k_0)}(X_i)}{k_0!}(d_i-\mu_i)^{k_0}\right].
\end{eqnarray}
Then
\begin{eqnarray}\nonumber
\frac{1}{n}\sum_{i}\mathbb{E}[t_i]\mathbb{E}\left[\frac{1}{d_i(d_i-1)}\right]&=&\frac{1}{n}\sum_{i}\frac{\mathbb{E}[t_i]}{\mu_i(\mu_i-1)}+\sum_{s=2}^{k_0-1}\frac{1}{n}\sum_{i}\mathbb{E}[t_i]\frac{h^{(s)}(\mu_i)}{s!}\mathbb{E}[(d_i-\mu_i)^s]\\ \label{eedd}
&&+\frac{1}{n}\sum_{i}\mathbb{E}[t_i]\mathbb{E}\left[\frac{h^{(k_0)}(X_i)}{k_0!}(d_i-\mu_i)^{k_0}\right].
\end{eqnarray}

Note that
\begin{eqnarray}\nonumber
\frac{1}{n}\sum_{i}\frac{2\mu_i-1}{\mu_i^2(\mu_i-1)^2}\mathbb{E}[t_i(d_i-\mu_i)]&=&\frac{1}{n}\sum_{i\neq j\neq k}\frac{2\mu_i-1}{\mu_i^2(\mu_i-1)^2}\mathbb{E}[A_{ij}\bar{A}_{ij}A_{jk}A_{ki}]\\ \nonumber
&&+\frac{1}{n}\sum_{i\neq j\neq k}\frac{2\mu_i-1}{\mu_i^2(\mu_i-1)^2}\mathbb{E}[A_{ij}A_{jk}A_{ki}\bar{A}_{ik}]\\ \label{teq6}
&=&\frac{2}{n}\sum_{i\neq j\neq k}\frac{2\mu_i-1}{\mu_i^2(\mu_i-1)^2}\mu_{ij}(1-\mu_{ij})\mu_{ik}\mu_{jk}.
\end{eqnarray}

Let $\delta_n=\left(\log(np_n)\right)^{-2}$. Then
\begin{eqnarray}\nonumber
\mathbb{E}\left[\big|h^{(s)}(X_i)(d_i-\mu_i)^{s}\big|\right]
&=&\mathbb{E}\left[\big|h^{(s)}(X_i)(d_i-\mu_i)^{s}\big|I[X_i\leq \delta_nnp_n]\right]\\ \label{eq9}
&&+\mathbb{E}\left[\big|h^{(s)}(X_i)(d_i-\mu_i)^{s}\big|I[X_i> \delta_nnp_n]\right].
\end{eqnarray}
For $x> \delta_nnp_n$, it is easy to verify that
\begin{eqnarray}\label{neq7}
|h^{(k)}(x)|&\leq&k!\frac{\sum_{t=1}^{k+1}\binom{k+1}{t}x^{k+1-t}}{(x-1)^{k+1}x^{k+1}}=O\left(\frac{1}{\big(\delta_nnp_n\big)^{k+2}}\right).
\end{eqnarray}
Then by Lemma \ref{lem2}, we get
\begin{eqnarray}\nonumber
\mathbb{E}\left[\big|h^{(s)}(X_i)(d_i-\mu_i)^{s}\big|I[X_i> \delta_nnp_n]\right]&=&O\left(\frac{\sqrt{\mathbb{E}\left[(d_i-\mu_i)^{2s}\right]}}{\big(\delta_nnp_n\big)^{s+2}}\right)\\ \label{eq7}
&=&O\left(\frac{1}{\delta_n^{s+2}\big(np_n\big)^{\frac{s}{2}+2}}\right).
\end{eqnarray}
If $X_i\leq \delta_nnp_n$, then $d_i\leq \delta_nnp_n$. Otherwise, $X_i$ can not be between $d_i$ and $\mu_i$. In this case, $|h^{(k)}(x)|=O(1)$. By Lemma \ref{lem1}, one has 
\begin{eqnarray}\nonumber
&&\mathbb{E}\left[\big|h^{(s)}(X_i)(d_i-\mu_i)^{s}\big|I[X_i\leq \delta_nnp_n]\right]\\ \nonumber
&\leq&O(1)\mathbb{E}\left[\big|(d_i-\mu_i)^{s}\big|I[d_i\leq \delta_nnp_n]\right]\\ \nonumber
&\leq&O(n^{s}))\sum_{t=2}^{\delta_nnp_n}\mathbb{P}(d_i=t)\\ \label{eq8}
&=&e^{-np_n\beta(1+o(1))}.
\end{eqnarray}
Note that $k_0\geq4$. Hence
\begin{eqnarray}\label{teq5}
\left|\frac{1}{n}\sum_{i}\mathbb{E}[t_i]\mathbb{E}\left[\frac{h^{(k_0)}(X_i)}{k_0!}(d_i-\mu_i)^{k_0}\right]\right|=O\left(\frac{p_n}{\delta_n^{k_0+2}\big(np_n\big)^{\frac{k_0}{2}}}\right)=o\left(\frac{1}{n\sqrt{np_n}}+\frac{\sqrt{p_n}}{n}\right).
\end{eqnarray}
Similarly, we have
\begin{eqnarray*} 
\left|\frac{1}{n}\sum_{i}\mathbb{E}\left[t_i\frac{h^{(k_0)}(X_i)}{k_0!}(d_i-\mu_i)^{k_0}\right]\right|&\leq&\frac{1}{n}\sum_{i}\frac{1}{k_0!}\sqrt{\mathbb{E}[t_i^2]\mathbb{E}\left[(h^{(k_0)}(X_i))^2(d_i-\mu_i)^{2k_0}\right]}\\
&=&o\left(\frac{1}{n\sqrt{np_n}}+\frac{\sqrt{p_n}}{n}\right).
\end{eqnarray*}

For $s\geq4$, 
\begin{eqnarray}\nonumber
\frac{1}{n}\sum_{i}\frac{|h^{(s)}(\mu_i)|}{s!}\mathbb{E}[|t_i(d_i-\mu_i)^s|]&\leq& \frac{1}{n}\sum_{i}\frac{|h^{(s)}(\mu_i)|}{s!}\sqrt{\mathbb{E}[t_i^2]\mathbb{E}[(d_i-\mu_i)^{2s}]}\\ \nonumber
&=&O\left(\frac{p_n}{(np_n)^{\frac{s}{2}}}\right)\\ \label{teq2}
&=&o\left(\frac{1}{n\sqrt{np_n}}+\frac{\sqrt{p_n}}{n}\right).
\end{eqnarray}

For $s=2$, we have
\begin{eqnarray}\nonumber
\frac{1}{2n}\sum_{i}h^{(2)}(\mu_i)\mathbb{E}[t_i(d_i-\mu_i)^2]&=&\frac{1}{2n}\sum_{i\neq j\neq k}h^{(2)}(\mu_i)\mathbb{E}[A_{ij}A_{ik}A_{jk}\bar{A}_{ij}\bar{A}_{ik}]\\ \nonumber
&&+\frac{1}{2n}\sum_{i\neq j\neq k\neq s}h^{(2)}(\mu_i)\mathbb{E}[A_{ij}A_{ik}A_{jk}\bar{A}_{is}^2]\\ \nonumber
&&+\frac{1}{2n}\sum_{i\neq j\neq k}h^{(2)}(\mu_i)\mathbb{E}[A_{ij}A_{ik}A_{jk}\bar{A}_{ij}^2]\\ \nonumber
&&+\frac{1}{2n}\sum_{i\neq j\neq k}h^{(2)}(\mu_i)\mathbb{E}[A_{ij}A_{ik}A_{jk}\bar{A}_{ik}^2]\\ \nonumber
&=&O\left(\frac{1}{n(np_n)}\right)+\frac{1}{2n}\sum_{i\neq j\neq k\neq s}h^{(2)}(\mu_i)\mu_{ij}\mu_{ik}\mu_{jk}\mu_{is}(1-\mu_{is}).
\end{eqnarray}
Moreover,
\begin{eqnarray}\nonumber
&&\frac{1}{2n}\sum_{i}h^{(2)}(\mu_i)\mathbb{E}[t_i]\mathbb{E}[(d_i-\mu_i)^2]\\ \nonumber
&=& \frac{1}{2n}\sum_{i\neq j\neq k\neq s}h^{(2)}(\mu_i)\mu_{ij}\mu_{ik}\mu_{jk}\mu_{is}(1-\mu_{is})+\frac{1}{2n}\sum_{i\neq j\neq k}h^{(2)}(\mu_i)\mu_{ij}\mu_{ik}\mu_{jk}\mu_{ij}(1-\mu_{ij})\\ \nonumber
&&+\frac{1}{2n}\sum_{i\neq j\neq k}h^{(2)}(\mu_i)\mu_{ij}\mu_{ik}\mu_{jk}\mu_{ik}(1-\mu_{ik})\\ \nonumber
&=& \frac{1}{2n}\sum_{i\neq j\neq k\neq s}h^{(2)}(\mu_i)\mu_{ij}\mu_{ik}\mu_{jk}\mu_{is}(1-\mu_{is})+O\left(\frac{p_n}{n(np_n)}\right).
\end{eqnarray}
Hence,
\begin{eqnarray}\label{teq1}
\frac{1}{2n}\sum_{i}h^{(2)}(\mu_i)\mathbb{E}[t_i(d_i-\mu_i)^2]-\frac{1}{2n}\sum_{i}h^{(2)}(\mu_i)\mathbb{E}[t_i]\mathbb{E}[(d_i-\mu_i)^2]
&=&O\left(\frac{1}{n(np_n)}\right).
\end{eqnarray}

Let $s=3$. In this case, $|\mathbb{E}[(d_i-\mu_i)^3]|=O(np_n)$. Then
\begin{eqnarray*}\label{teq3}
\frac{1}{6n}\sum_{i}\big|h^{(3)}(\mu_i)\mathbb{E}[t_i]\mathbb{E}[(d_i-\mu_i)^3]\big|=O\left(\frac{1}{n(np_n)}\right).
\end{eqnarray*}
Moreover, 
\begin{eqnarray*}\nonumber
\frac{1}{6n}\sum_{i}h^{(3)}(\mu_i)\mathbb{E}[t_i(d_i-\mu_i)^3]
&=& \frac{1}{6n}\sum_{i\neq j\neq k,r\neq s\neq t}h^{(3)}(\mu_i)\mathbb{E}[A_{ij}A_{ik}A_{jk}\bar{A}_{ir}\bar{A}_{is}\bar{A}_{it}]\\ \nonumber
&&+\frac{1}{6n}\sum_{i\neq j\neq k\neq s}h^{(3)}(\mu_i)\mathbb{E}[A_{ij}A_{ik}A_{jk}\bar{A}_{is}^2\bar{A}_{ij}]\\ \nonumber
&&+\frac{1}{6n}\sum_{i\neq j\neq k, s}h^{(3)}(\mu_i)\mathbb{E}[A_{ij}A_{ik}A_{jk}\bar{A}_{is}^3]\\ \nonumber
&=&O\left(\frac{1}{n(np_n)}\right).
\end{eqnarray*}
Hence,
\begin{eqnarray}\label{teq4}
\frac{1}{6n}\sum_{i}h^{(3)}(\mu_i)\mathbb{E}[t_i(d_i-\mu_i)^3]-\frac{1}{6n}\sum_{i}h^{(3)}(\mu_i)\mathbb{E}[t_i]\mathbb{E}[(d_i-\mu_i)^3]
&=&O\left(\frac{1}{n(np_n)}\right).
\end{eqnarray}

By (\ref{00eq2}),  (\ref{eedd}), (\ref{teq6}), (\ref{teq5}), (\ref{teq2}), (\ref{teq1}) and (\ref{teq4}), the proof is complete.

\qed

\begin{Lemma}\label{nlem1}
Let $a_i=\mathbb{E}\left[\frac{1}{d_i(d_i-1)}\right]$, $a_{ijk}=a_i+a_j+a_k$ and 
\[\sigma_{1n}^2=\frac{4}{n^2}\sum_{i<j<k}a_{ijk}^2\mu_{ij}(1-\mu_{ij})\mu_{jk}(1-\mu_{jk})\mu_{ki}(1-\mu_{ki}).\]
\[\sigma_{2n}^2=\frac{4}{n^2}\sum_{i<j}e_{ij}^2\mu_{ij}(1-\mu_{ij}).\]
For any fixed constants $\lambda_1,\lambda_2$ with $\lambda_1^2+\lambda_2^2=1$, we have
    \[\frac{\lambda_1\frac{2}{n}\sum_{i<j<k}a_{ijk}\bar{A}_{ij}\bar{A}_{jk}\bar{A}_{ki}+\lambda_2\frac{2}{n}\sum_{i<j}e_{ij}\bar{A}_{ij}}{\sqrt{\lambda_1^2\sigma_{1n}^2+\lambda_2^2\sigma_{2n}^2}}\Rightarrow \mathcal{N}(0,1).\]
\end{Lemma}

To prove  Lemma \ref{nlem1}, we need the following proposition.
    
\begin{Proposition}[\cite{HH14}]\label{martingale}
 Suppose that for every $n\in\mathbb{N}$ and $k_n\rightarrow\infty$ the random variables $X_{n,1},\dots,X_{n,k_n}$ are a martingale difference sequence relative to an arbitrary filtration $\mathcal{F}_{n,1}\subset\mathcal{F}_{n,2}$ $\subset$ $\dots$ $\subset\mathcal{F}_{n,k_n}$. If (I) $\sum_{i=1}^{k_n}\mathbb{E}(X_{n,i}^2|\mathcal{F}_{n,i-1})\rightarrow 1$ in probability,
 (II) $\sum_{i=1}^{k_n}\mathbb{E}(X_{n,i}^2I[|X_{n,i}|>\epsilon]|\mathcal{F}_{n,i-1})\rightarrow 0$ in probability for every $\epsilon>0$,
\noindent then $\sum_{i=1}^{k_n}X_{n,i}\rightarrow N(0,1)$ in distribution.
\end{Proposition}

\medskip

\noindent
{\bf Proof of Lemma \ref{nlem1}:} We employ Proposition \ref{martingale} to prove Lemma \ref{nlem1}. Let $\sigma_n^2=\lambda_1^2\sigma_{1n}^2+\lambda_2^2\sigma_{2n}^2$,
\[Y_t=\frac{\lambda_1\frac{2}{n}\sum_{1\leq i<j<k\leq t}(a_{ijk})\bar{A}_{ij}\bar{A}_{jk}\bar{A}_{ki}+\lambda_2\frac{2}{n}\sum_{i<j\leq t}e_{ij}\bar{A}_{ij}}{\sigma_n},\]
for $3\leq t\leq n$, and $Y_2=0$. Then $\{Y_t\}_{t=2}^n$ is a martingale and 
\[Y_n=\frac{\lambda_1\frac{2}{n}\sum_{1\leq i<j< k}a_{ijk}\bar{A}_{ij}\bar{A}_{jk}\bar{A}_{ki}+\lambda_2\frac{2}{n}\sum_{i<j}e_{ij}\bar{A}_{ij}}{\sigma_n}.\]

Let $X_t=Y_t-Y_{t-1}$ and $F_t=\{A_{ij},1\leq i<j\leq t\}$ for $t\geq3$. It is easy to verify that
\[X_t=\frac{\lambda_1\frac{2}{n}\sum_{1\leq i<j< t}a_{ijt}\bar{A}_{ij}\bar{A}_{jt}\bar{A}_{ti}+\lambda_2\frac{2}{n}\sum_{i<t}e_{it}\bar{A}_{it}}{\sigma_n},\]
and
 $\mathbb{E}[X_{t}|F_{t-1}]=0$. Hence, $\{X_t\}_{t=3}^n$ is a martingale difference. 
 
 Now we verify the two conditions in Proposition \ref{martingale}. Firstly, we prove condition $(I)$ is satisfied.
It is easy to get $\mathbb{E}\big(Y_{t}|F_{t-1}\big)=Y_{t-1}$. By the property of conditional expectation, one has
\begin{eqnarray*}
\mathbb{E}\left[\sum_{t=3}^n\mathbb{E}\Big(X^2_{t}|F_{t-1}\Big)\right]&=&\mathbb{E}\left[\sum_{t=3}^n\mathbb{E}\Big((Y^2_{t}-2Y_{t}Y_{t-1}+Y^2_{t-1})|F_{t-1}\Big)\right]\\
&=&\mathbb{E}\left[\sum_{t=3}^n\mathbb{E}\big[(Y^2_{t}-Y^2_{t-1})|F_{t-1}\big]\right]=\mathbb{E}[Y^2_{n}]=1.
\end{eqnarray*}
To prove condition $(I)$ holds, we only need to show that \[\mathbb{E}\Big(\sum_{t=3}^n\mathbb{E}[X_{t}^2|F_{t-1}]\Big)^2=1+o(1).\] 

For convenience, let $\sigma_{ij}^2=\mu_{ij}(1-\mu_{ij})$.
Given $t\in\{3,4,\dots,n\}$, one has
\begin{eqnarray*}
\mathbb{E}[X_{t}^2|F_{t-1}]
&=&\frac{4\lambda_1^2}{n^2\sigma_n^2}\sum_{\substack{1\leq i_1<j_1<t\\1\leq i_2<j_2<t}}a_{i_1j_1t}a_{i_2j_2t}\mathbb{E}\Big[\bar{A}_{i_1j_1}\bar{A}_{j_1t}\bar{A}_{ti_1}\bar{A}_{i_2j_2}\bar{A}_{j_2t}\bar{A}_{ti_2}|F_{t-1}\Big]\\
&&+\frac{4\lambda_2^2}{n^2\sigma_n^2}\sum_{\substack{1\leq i_1<t\\1\leq i_2<t}}e_{i_1t}e_{i_2t}\mathbb{E}\Big[\bar{A}_{i_1t}\bar{A}_{i_2t}|F_{t-1}\Big]    \\
&&+  \frac{8\lambda_1\lambda_2}{n^2\sigma_n^2} \sum_{\substack{1\leq i_1<j_1<t\\1\leq i_2<t}}a_{i_1j_1t}e_{i_2t}\mathbb{E}\Big[\bar{A}_{i_1j_1}\bar{A}_{j_1t}\bar{A}_{ti_1}\bar{A}_{i_2t}|F_{t-1}\Big]         \\
&=&\frac{4\lambda_1^2}{n^2\sigma_n^2}\sum_{\substack{1\leq i<j<t\\}}a_{ijt}^2\bar{A}_{ij}^2\sigma_{jt}^2\sigma_{it}^2+\frac{4\lambda_2^2}{n^2\sigma_n^2}\sum_{\substack{1\leq i<t}}e_{it}^2\sigma_{it}^2
\end{eqnarray*}
Then we have
\begin{eqnarray*}
\sum_{t=3}^n\mathbb{E}[X_{t}^2|F_{t-1}]
&=&\frac{4\lambda_1^2}{n^2\sigma_n^2}\sum_{\substack{1\leq i<j<t\leq n\\}}a_{ijt}^2\bar{A}_{ij}^2\sigma_{jt}^2\sigma_{it}^2+\frac{\lambda_2^2\sigma_{2n}^2}{\sigma_n^2}.
\end{eqnarray*}
Note that $A_{ij}$ $(1\leq i<j\leq n)$ are independent,
\[\mathbb{E}\Big[\bar{A}_{ij}^2(A_{kl}-\mu_{kl})^2\Big]=O(p_n),\hskip  1cm if\ \{i,j\}=\{k,l\},\]
and 
\[\mathbb{E}\Big[\bar{A}_{ij}^2(A_{kl}-\mu_{kl})^2\Big]=\sigma_{ij}^2\sigma_{kl}^2, \hskip 1cm if\ |\{i,j\}\cap\{k,l\}|\leq1.\]
 Then we have
\begin{eqnarray}\nonumber
&&\mathbb{E}\left[\Big(\frac{4\lambda_1^2}{n^2\sigma_n^2}\sum_{\substack{1\leq i<j<t\\}}a_{ijt}^2\bar{A}_{ij}^2\sigma_{jt}^2\sigma_{it}^2\Big)^2\right]\\
&=&\frac{16\lambda_1^4}{n^4\sigma_n^4}\sum_{\substack{1\leq i<j<t\leq n\\1\leq k<l<s\leq n}}a_{ijt}^2a_{kls}^2\sigma_{jt}^2\sigma_{it}^2\sigma_{ks}^2\sigma_{ls}^2\mathbb{E}[\bar{A}_{ij}^2\bar{A}_{kl}^2]\\  \nonumber
&=&\frac{16\lambda_1^4}{n^4\sigma_n^4}\sum_{\substack{1\leq i<j<t\leq n\\1\leq k<l<s\leq n\\|\{i,j\}\cap\{k,l\}|\leq1}}a_{ijt}^2a_{kls}^2\sigma_{jt}^2\sigma_{it}^2\sigma_{ks}^2\sigma_{ls}^2\mathbb{E}[\bar{A}_{ij}^2\bar{A}_{kl}^2]+\frac{16\lambda_1^4}{n^4\sigma_n^4}\sum_{\substack{1\leq i<j<t\leq n\\1\leq k<l<s\leq n\\ \{i,j\}=\{k,l\}}}a_{ijt}^2a_{kls}^2\sigma_{jt}^2\sigma_{it}^2\sigma_{ks}^2\sigma_{ls}^2\mathbb{E}[\bar{A}_{ij}^4]\\ \nonumber
&=& \frac{16\lambda_1^4}{n^4\sigma_n^4}\sum_{\substack{1\leq i<j<t\leq n\\1\leq k<l<s\leq n\\|\{i,j\}\cap\{k,l\}|\leq1}}a_{ijt}^2a_{klt}^2\sigma_{ij}^2\sigma_{jt}^2\sigma_{it}^2\sigma_{kl}^2\sigma_{ks}^2\sigma_{ls}^2+\frac{16\lambda_1^4}{n^4\sigma_n^4}\sum_{\substack{1\leq i<j<t\leq n\\1\leq k<l<s\leq n\\ \{i,j\}=\{k,l\}}}a_{ijt}^2a_{klt}^2\mathbb{E}[\bar{A}_{ij}^4]\sigma_{jt}^2\sigma_{it}^2\sigma_{ks}^2\sigma_{ls}^2\\  \nonumber
&=&\frac{16\lambda_1^4}{n^4\sigma_n^4}\sum_{\substack{1\leq i<j<t\leq n\\1\leq k<l<s\leq n\\|\{i,j\}\cap\{k,l\}|\leq1}}a_{ijt}^2a_{klt}^2\sigma_{ij}^2\sigma_{jt}^2\sigma_{it}^2\sigma_{kl}^2\sigma_{ks}^2\sigma_{ls}^2+\frac{16\lambda_1^4}{n^4\sigma_n^4}\sum_{\substack{1\leq i<j<t\leq n\\1\leq k<l<s\leq n\\\{i,j\}=\{k,l\}}}a_{ijt}^2a_{klt}^2\sigma_{ij}^2\sigma_{jt}^2\sigma_{it}^2\sigma_{kl}^2\sigma_{ks}^2\sigma_{ls}^2\\ \nonumber
&&+\frac{16\lambda_1^4}{n^4\sigma_n^4}\sum_{\substack{1\leq i<j<t\leq n\\1\leq k<l<s\leq n\\ \{i,j\}=\{k,l\}}}a_{ijt}^2a_{klt}^2\Big(\mathbb{E}[\bar{A}_{ij}^4]\sigma_{jt}^2\sigma_{it}^2\sigma_{ks}^2\sigma_{ls}^2-\sigma_{ij}^2\sigma_{jt}^2\sigma_{it}^2\sigma_{kl}^2\sigma_{ks}^2\sigma_{ls}^2)\Big)\\
&=&\frac{\lambda_1^4\sigma_{1n}^4}{\sigma_n^4}+O\Big(\frac{1}{(np_n)^4}\Big).
\end{eqnarray}

 Then
\begin{eqnarray}\nonumber
\mathbb{E}\Big(\sum_{t=1}^n\mathbb{E}[X_{n,t}^2|F_{t-1}]\Big)^2&=&\mathbb{E}\left[\Big(\frac{4\lambda_1^2}{n^2\sigma_n^2}\sum_{\substack{1\leq i<j<t\\}}a_{ijt}^2\bar{A}_{ij}^2\sigma_{jt}^2\sigma_{it}^2\Big)^2\right]+\left(\frac{\lambda_2^2\sigma_{2n}^2}{\sigma_n^2}\right)^2\\  \nonumber
&&+2\left(\frac{\lambda_2^2\sigma_{2n}^2}{\sigma_n^2}\right)\frac{4\lambda_1^2}{n^2\sigma_n^2}\sum_{\substack{1\leq i<j<t\leq n\\}}a_{ijt}^2\mathbb{E}[\bar{A}_{ij}^2]\sigma_{jt}^2\sigma_{it}^2\\ \nonumber
&=&\frac{\lambda_1^4\sigma_{1n}^4}{\sigma_n^4}+\left(\frac{\lambda_2^2\sigma_{2n}^2}{\sigma_n^2}\right)^2+2\left(\frac{\lambda_1^2\sigma_{1n}^2}{\sigma_n^2}\right)\left(\frac{\lambda_2^2\sigma_{2n}^2}{\sigma_n^2}\right)+o(1)\\ \nonumber
&=&1+o(1).
\end{eqnarray}

 Now we check condition $(II)$ in Proposition \ref{martingale}. 
Let $\epsilon$ be a fixed positive constant. By  the Cauchy-Schwarz inequality and Markov's inequality, we have
\begin{eqnarray}\nonumber
&&\mathbb{E}\Big[\sum_{t=3}^n\mathbb{E}\big[X_{t}^2I[|X_{t}|>\epsilon|F_{t-1}\big]\Big]\\    \nonumber
&\leq&\mathbb{E}\Big[\sum_{t=3}^n\sqrt{\mathbb{E}\big[X_{t}^4|F_{t-1}\big]\mathbb{P}[|X_{t}|>\epsilon|F_{t-1}\big]}\Big]\\  \nonumber
&\leq&\mathbb{E}\Bigg[\frac{4}{\epsilon^2n^4\sigma_n^4}\sum_{t=1}^n\mathbb{E}\Big[\Big(\lambda_1\sum_{1\leq i<j< t}a_{ijt}\bar{A}_{ij}\bar{A}_{jt}\bar{A}_{ti}+\lambda_2\sum_{i<t}e_{it}\bar{A}_{it}\Big)^4\Big|F_{t-1}\Big]\Bigg]\\   \nonumber
&\leq&\frac{32}{\epsilon^2n^4\sigma_n^4}\sum_{t=3}^n\sum_{\substack{1\leq i_1<j_1<t\\ 1\leq i_2<j_2<t\\ 1\leq i_3<j_3<t\\ 1\leq i_4<j_4<t}}\lambda_1^4a_{i_1j_1t}a_{i_2j_2t}a_{i_3j_3t}a_{i_4j_4t}\mathbb{E}\Big[\bar{A}_{i_1j_1}\bar{A}_{j_1t}\bar{A}_{ti_1}\bar{A}_{i_2j_2}\bar{A}_{j_2t}\bar{A}_{ti_2}\\
&&\times\bar{A}_{i_3j_3}\bar{A}_{j_3t}\bar{A}_{ti_3}\bar{A}_{i_4j_4}\bar{A}_{j_4t}\bar{A}_{ti_4}\Big]\\   \nonumber
&&+\frac{32}{\epsilon^2n^4\sigma_n^4}\sum_{t=3}^n\sum_{\substack{1\leq i_1<t\\ 1\leq i_2<t\\ 1\leq i_3<t\\ 1\leq i_4<t}}\lambda_2^4e_{i_1t}e_{i_2t}e_{i_3t}e_{i_4t}\mathbb{E}\Big[\bar{A}_{i_1t}\bar{A}_{i_2t}\bar{A}_{i_3t}\bar{A}_{i_4t}\Big]\\
&\leq&\frac{32C}{\epsilon^2n^4\sigma_n^4}\sum_{t=3}^n\sum_{\substack{1\leq i_1<j_1<t\\ 1\leq i_2<j_2<t}}a_{i_1j_1t}^2a_{i_2j_2t}^2\mathbb{E}\big[\bar{A}_{i_1j_1}^2\bar{A}_{j_1t}^2\bar{A}_{i_1t}^2 \bar{A}_{i_2j_2}^2\bar{A}_{j_2t}^2\bar{A}_{ti_2}^2\big]\\   \nonumber
&&+\frac{32C}{\epsilon^2n^4\sigma_n^4}\sum_{t=3}^n\sum_{\substack{1\leq i_1<j_1<t}}a_{i_1j_1t}^4\mathbb{E}\big[\bar{A}_{i_1j_1}^4\bar{A}_{j_1t}^4\bar{A}_{i_1t}^4\big] \\   \nonumber
&&+\frac{32}{\epsilon^2n^4\sigma_n^4}\sum_{t=3}^n\sum_{\substack{1\leq i_1<t\\ 1\leq i_2<t}}\lambda_2^4e_{i_1t}^2e_{i_2t}^2\mathbb{E}\Big[\bar{A}_{i_1t}^2\bar{A}_{i_2t}^2\Big]+\frac{32}{\epsilon^2n^4\sigma_n^4}\sum_{t=3}^n\sum_{\substack{1\leq i_1<t}}\lambda_2^4e_{i_1t}^4\mathbb{E}\Big[\bar{A}_{i_1t}^4\Big]\\ \nonumber
&=&O\Big(\frac{n^5p_n^6}{\epsilon^2(np_n)^8n^4\sigma_n^4}\Big)+O\Big(\frac{n^3p_n^3}{\epsilon^2(np_n)^8n^4\sigma_n^4}\Big)+O\Big(\frac{p_n^2}{\epsilon^2(np_n)^8n^5\sigma_n^4}\Big)\\ \nonumber
&=&o(1).
\end{eqnarray}
 Then the desired result follows from Proposition \ref{martingale}.

\qed

\subsection{Proof of Theorem \ref{mainthm}}
By Lemma \ref{lem3}, straightforward calculation yields
\begin{eqnarray}\label{eq1}
\overline{\mathcal{C}}_n-\mathbb{E}[\overline{\mathcal{C}}_n]
=Y_{1n}+Y_{2n}+Y_{3n}+o\left(\frac{1}{n\sqrt{np_n}}+\frac{\sqrt{p_n}}{n}\right),
\end{eqnarray}
where
\begin{eqnarray}\nonumber
Y_{1n}
&=&\frac{1}{n}\sum_{i=1}^n(t_i-\mathbb{E}[t_i])\mathbb{E}\left[\frac{1}{d_i(d_i-1)}\right],\\ \nonumber
Y_{2n}&=&\frac{1}{n}\sum_{i=1}^n(t_i-\mathbb{E}[t_i])\left(\frac{1}{d_i(d_i-1)}-\mathbb{E}\left[\frac{1}{d_i(d_i-1)}\right]\right)\\
&&+\frac{2}{n}\sum_{i\neq j\neq k}\frac{2\mu_i-1}{\mu_i^2(\mu_i-1)^2}\mathbb{E}[\bar{A}_{ij}^2]\mu_{ik}\mu_{jk},\\
Y_{3n}&=&\frac{1}{n}\sum_{i=1}^n\mathbb{E}[t_i]\left(\frac{1}{d_i(d_i-1)}-\mathbb{E}\left[\frac{1}{d_i(d_i-1)}\right]\right).
\end{eqnarray}

Next we find the order of
$Y_{1n}$, $Y_{2n}$ and $Y_{3n}$.

\medskip

{\bf (a) Order of $Y_{3n}$.} Firstly, we find the order of $Y_{3n}$. 
By (\ref{0eq2}) and (\ref{edd}), one has
\begin{eqnarray}\nonumber
&&\frac{1}{d_i(d_i-1)}-\mathbb{E}\left[\frac{1}{d_i(d_i-1)}\right]\\ \nonumber
&=&-\frac{2\mu_i-1}{\mu_i^2(\mu_i-1)^2}(d_i-\mu_i)+\sum_{s=2}^{k_0-1}\frac{h^{(s)}(\mu_i)}{s!}\Big[(d_i-\mu_i)^s-\mathbb{E}[(d_i-\mu_i)^s]\Big]\\ \label{eq3}
&&+\frac{h^{(k_0)}(X_i)}{k_0!}(d_i-\mu_i)^{k_0}-\mathbb{E}\left[\frac{h^{(k_0)}(X_i)}{k_0!}(d_i-\mu_i)^{k_0}\right].
\end{eqnarray}
Then we can express $Y_{3n}$ as follows.
\begin{eqnarray}\nonumber
Y_{3n}&=&-\frac{1}{n}\sum_{i=1}^n\frac{\mathbb{E}[t_i](2\mu_i-1)}{\mu_i^2(\mu_i-1)^2}(d_i-\mu_i)+\sum_{s=2}^{k_0-1}\frac{1}{n}\sum_{i=1}^n\frac{h^{(s)}(\mu_i)\mathbb{E}[t_i]}{s!}\Big[(d_i-\mu_i)^s-\mathbb{E}[(d_i-\mu_i)^s]\Big]\\ \label{eq6}
&&+\frac{1}{n}\sum_{i=1}^n\frac{h^{(k_0)}(X_i)\mathbb{E}[t_i]}{k_0!}(d_i-\mu_i)^{k_0}-\frac{1}{n}\sum_{i=1}^n\mathbb{E}[t_i]\mathbb{E}\left[\frac{h^{(k_0)}(X_i)}{k_0!}(d_i-\mu_i)^{k_0}\right].
\end{eqnarray}

Now we find the order of each term in (\ref{eq6}) of $Y_{3n}$. Recall that $A_{ij} (1\leq i<j\leq n)$ are independent and $\mathbb{E}[A_{ij}]=\mu_{ij}$. Hence, for indices $i,j (i\neq j)$ and $k,l (k\neq l)$, we have
\begin{equation}\label{eq5}
 \mathbb{E}[\bar{A}_{ij}\bar{A}_{kl}] =
    \begin{cases}
      \mu_{ij}(1-\mu_{ij}) & \text{if $\{i,j\}=\{k,l\}$}, \\
      0 & \text{if $\{i,j\}\neq\{k,l\}$}.     
    \end{cases}       
\end{equation}

The first term in (\ref{eq6}) is equal to
\begin{eqnarray}\label{eq55}
-\frac{1}{n}\sum_{i=1}^n\frac{\mathbb{E}[t_i](2\mu_i-1)}{\mu_i^2(\mu_i-1)^2}(d_i-\mu_i)=-\frac{1}{n}\sum_{i\neq j}^n\frac{\mathbb{E}[t_i](2\mu_i-1)}{\mu_i^2(\mu_i-1)^2}\bar{A}_{ij}.
\end{eqnarray}
By (\ref{eq5}) and (\ref{eq55}), we have
\begin{eqnarray}\nonumber
&&\mathbb{E}\left[\left(\frac{1}{n}\sum_{i=1}^n\frac{\mathbb{E}[t_i](2\mu_i-1)}{\mu_i^2(\mu_i-1)^2}(d_i-\mu_i)\right)^2\right]\\  \nonumber
&=&\frac{1}{n^2}\sum_{i\neq j,s\neq t}^n\frac{\mathbb{E}[t_i](2\mu_i-1)}{\mu_i^2(\mu_i-1)^2}\frac{\mathbb{E}[t_s](2\mu_s-1)}{\mu_s^2(\mu_s-1)^2}\mathbb{E}\left[\bar{A}_{ij}\bar{A}_{st}\right]\\ \nonumber
&=&\frac{1}{n^2}\sum_{i\neq j}^n\frac{(\mathbb{E}[t_i](2\mu_i-1))^2}{\mu_i^4(\mu_i-1)^4}\mathbb{E}\left[\bar{A}_{ij}^2\right]\\ \nonumber
&&+\frac{1}{n^2}\sum_{i\neq j}^n\frac{\mathbb{E}[t_i](2\mu_i-1)}{\mu_i^2(\mu_i-1)^2}\frac{\mathbb{E}[t_j](2\mu_j-1)}{\mu_j^2(\mu_j-1)^2}\mathbb{E}\left[\bar{A}_{ij}^2\right]\\ \nonumber
&=&\frac{1}{n^2}\sum_{i\neq j}^n\frac{(\mathbb{E}[t_i](2\mu_i-1))^2}{\mu_i^4(\mu_i-1)^4}\mu_{ij}(1-\mu_{ij})+\frac{1}{n^2}\sum_{i\neq j}^n\frac{\mathbb{E}[t_i](2\mu_i-1)}{\mu_i^2(\mu_i-1)^2}\frac{\mathbb{E}[t_j](2\mu_j-1)}{\mu_j^2(\mu_j-1)^2}\mu_{ij}(1-\mu_{ij})\\ \label{eq4}
&=&\Theta\left(\frac{p_n}{n^2}\right).
\end{eqnarray}
Here the last equality follows from the facts that $\beta p_n\leq\mu_{ij}\leq p_n$, $\beta np_n\leq\mu_i\leq np_n$ and $\mathbb{E}[t_i]=\Theta(n^2p_n^3)$. 
Then the first term of $Y_{3n}$ is of order $\frac{\sqrt{p_n}}{n}$.

Now we show the last two terms of $Y_{3n}$ are $o_P\left(\frac{\sqrt{p_n}}{n}+\frac{1}{n\sqrt{np_n}}\right)$.
Recall that $k_0=\lceil 1+\frac{1}{1-\alpha}\rceil +2$ in the proof of Lemma \ref{lem3}. Then $k_0>\max\{1+\frac{1}{1-\alpha},3\}$. By (\ref{eq9}), (\ref{eq7}) and (\ref{eq8}), we get
\begin{eqnarray}\nonumber
&&\left|\frac{1}{n}\sum_{i=1}^n\mathbb{E}[t_i]\left(\frac{h^{(k_0)}(X_i)}{k_0!}(d_i-\mu_i)^{k_0}-\mathbb{E}\left[\frac{h^{(k_0)}(X_i)}{k_0!}(d_i-\mu_i)^{k_0}\right]\right)\right|\\ \nonumber
&=&O_P\left(\frac{n^2p_n^3}{\delta_n^{k_0+2}\big(np_n\big)^{\frac{k_0}{2}+2}}+n^2p_n^3e^{-np_n\beta(1+o(1))}\right)\\ \label{eq10}
&=&o_P\left(\frac{\sqrt{p_n}}{n}+\frac{1}{n\sqrt{np_n}}\right).
\end{eqnarray}
Hence the last two terms of $Y_{3n}$ are $o_P\left(\frac{\sqrt{p_n}}{n}+\frac{1}{n\sqrt{np_n}}\right)$.

Now we show the second term of $Y_{3n}$ is $o_P\left(\frac{\sqrt{p_n}}{n}+\frac{1}{n\sqrt{np_n}}\right)$. In the expression of $(d_i-\mu_i)^s$ in (\ref{dseq}), suppose $\lambda_{t1}=\dots=\lambda_{tt_1}=1$ and $\lambda_{t,t_1+1}\geq2$, \dots, $\lambda_{tt}\geq2$ for $t_1\geq1$. If $s=t$, then $t_1=s\geq2$. If $s>t$, then $2(s-t)\geq2$. Hence, we have
\begin{eqnarray*}
&&\mathbb{E}\left[\left(\frac{1}{n}\sum_{j_1\neq j_2\neq \dots,j_t\neq i}\mathbb{E}[t_i]h^{(s)}(\mu_i)\prod_{l=1}^{t_1}\bar{A}_{ij_l}\prod_{l=t_1+1}^t\bar{A}_{ij_l}^{\lambda_{tl}}\right)^2\right]\\
&\leq&\frac{1}{n^2}\sum_{\substack{j_1\neq j_2\neq \dots,j_t\neq i\\ k_{t_1+1}\neq\dots\neq k_s\neq i }}\left(\mathbb{E}[t_i]h^{(s)}(\mu_i)\right)^2\mathbb{E}\left[\prod_{l=1}^{t_1}\bar{A}_{ij_l}^2\prod_{l=t_1+1}^t\bar{A}_{ij_l}^{\lambda_{tl}} \prod_{l=t_1+1}^s\bar{A}_{ik_l}^{\lambda_{tl}} \right]\\
&=&O\left(\frac{n^{2t-t_1+1}n^4p_n^6p_n^{2t-t_1}}{n^2(np_n)^{2s+4}}\right)\\
&=&O\left(\frac{p_n^3}{(np_n)^{2(s-t)+t_1+1}}\right)\\
&=&o\left(\frac{p_n}{n^2}+\frac{1}{n^2(np_n)}\right).
\end{eqnarray*}

Now suppose $\lambda_{tl}\geq2$ for all $l=1,2,\dots,t$ in (\ref{dseq}). In this case, $t\leq\frac{s}{2}$. Then $s-t\geq1$. Note that
\begin{eqnarray*}
\prod_{l=1}^t\bar{A}_{ij_l}^{\lambda_{tl}}-\prod_{l=1}^t\mathbb{E}[\bar{A}_{ij_l}^{\lambda_{tl}}]
&=&\sum_{t_1=1}^{t}C_{t_1}\prod_{l=1}^{t_1}\left(\bar{A}_{ij_l}^{\lambda_{tl}}-\mathbb{E}[\bar{A}_{ij_l}^{\lambda_{tl}}]\right)\prod_{l=t_1+1}^t\mathbb{E}[\bar{A}_{ij_l}^{\lambda_{tl}}],
\end{eqnarray*}
for some constants $C_{t_1}$.
Fix $t_1$ ($1\leq t_1\leq t$). Then
\begin{eqnarray*}
&&\mathbb{E}\left[\left(\frac{1}{n}\sum_{j_1\neq j_2\neq \dots,j_t\neq i}\mathbb{E}[t_i]h^{(s)}(\mu_i)\prod_{l=1}^{t_1}\left(\bar{A}_{ij_l}^{\lambda_{l}}-\mathbb{E}[\bar{A}_{ij_l}^{\lambda_{l}}]\right)\prod_{l=t_1+1}^t\mathbb{E}[\bar{A}_{ij_l}^{\lambda_{l}}]\right)^2\right]\\
&\leq&\frac{1}{n^2}\sum_{\substack{j_1\neq j_2\neq \dots,j_t\neq i\\ k_{t_1+1},k_t}}\left(\mathbb{E}[t_i]h^{(s)}(\mu_i)\right)^2\prod_{l=1}^{t_1}\left(\bar{A}_{ij_l}^{\lambda_{l}}-\mathbb{E}[\bar{A}_{ij_l}^{\lambda_{l}}]\right)^2\prod_{l=t_1+1}^t\mathbb{E}[\bar{A}_{ij_l}^{\lambda_{l}}]\\
&&\times \prod_{l=t_1+1}^t\mathbb{E}[\bar{A}_{ik_l}^{\lambda_{l}}]\\
&=&O\left(\frac{n^{2t-t_1+1}n^4p_n^6p_n^{2t-t_1}}{n^2(np_n)^{2s+4}}\right)\\
&=&O\left(\frac{p_n^3}{(np_n)^{2(s-t)+t_1+1}}\right)\\
&=&o\left(\frac{p_n}{n^2}+\frac{1}{n^2(np_n)}\right).
\end{eqnarray*}

Hence, the second term of $Y_{3n}$ is $o_P\left(\frac{\sqrt{p_n}}{n}+\frac{1}{n\sqrt{np_n}}\right)$. Then we get
\begin{eqnarray}\label{y3n}
Y_{3n}=-\frac{1}{n}\sum_{i=1}^n\frac{\mathbb{E}[t_i](2\mu_i-1)}{\mu_i^2(\mu_i-1)^2}(d_i-\mu_i)+o_P\left(\frac{\sqrt{p_n}}{n}+\frac{1}{n\sqrt{np_n}}\right).
\end{eqnarray}

{\bf (b) Order of  $Y_{1n}$.} Straightforward calculation yields
\begin{eqnarray}\nonumber
A_{ij}A_{jk}A_{ki}-\mu_{ij}\mu_{jk}\mu_{ki}
&=&\bar{A}_{ij}\bar{A}_{ik}\bar{A}_{jk}+\bar{A}_{ij}\bar{A}_{ik}\mu_{jk}+\bar{A}_{ij}\bar{A}_{jk}\mu_{ik}+\bar{A}_{ik}\bar{A}_{jk}\mu_{ij}\\ \label{neq15}
&&+\bar{A}_{jk}\mu_{ij}\mu_{ik}+\bar{A}_{ij}\mu_{jk}\mu_{ik}+\bar{A}_{ik}\mu_{jk}\mu_{ij}.
\end{eqnarray}
Then $Y_{1n}$ can be expressed as follows.
\begin{eqnarray}\nonumber
Y_{1n}
&=&\frac{1}{n}\sum_{i\neq j\neq k}\mathbb{E}\left[\frac{1}{d_i(d_i-1)}\right]\bar{A}_{ij}\bar{A}_{ik}\bar{A}_{jk}+\frac{1}{n}\sum_{i\neq j\neq k}\mathbb{E}\left[\frac{1}{d_i(d_i-1)}\right]\bar{A}_{ij}\bar{A}_{ik}\mu_{jk}\\ \nonumber
&&+\frac{1}{n}\sum_{i\neq j\neq k}\mathbb{E}\left[\frac{1}{d_i(d_i-1)}\right]\bar{A}_{ij}\bar{A}_{jk}\mu_{ik}+\frac{1}{n}\sum_{i\neq j\neq k}\mathbb{E}\left[\frac{1}{d_i(d_i-1)}\right]\bar{A}_{ik}\bar{A}_{jk}\mu_{ij}\\ \nonumber
&&+\frac{1}{n}\sum_{i\neq j\neq k}\mathbb{E}\left[\frac{1}{d_i(d_i-1)}\right]\bar{A}_{jk}\mu_{ij}\mu_{ik}+\frac{1}{n}\sum_{i\neq j\neq k}\mathbb{E}\left[\frac{1}{d_i(d_i-1)}\right]\bar{A}_{ij}\mu_{jk}\mu_{ik}\\ \label{neq10}
&&+\frac{1}{n}\sum_{i\neq j\neq k}\mathbb{E}\left[\frac{1}{d_i(d_i-1)}\right]\bar{A}_{ik}\mu_{jk}\mu_{ij}
\end{eqnarray}
We find the order of each term in $Y_{1n}$.
Recall that we denote $a_i=\mathbb{E}\left[\frac{1}{d_i(d_i-1)}\right]$. The second moment of the first term of (\ref{neq10}) is equal to 
\begin{eqnarray*}
&&\mathbb{E}\left[\left(\frac{1}{n}\sum_{i\neq j\neq k}a_i\bar{A}_{ij}\bar{A}_{ik}\bar{A}_{jk}\right)^2\right]\\
&=&\frac{1}{n^2}\sum_{i\neq j\neq k}(a_i^2+a_ia_j+a_ia_k)\mathbb{E}[\bar{A}_{ij}^2]\mathbb{E}[\bar{A}_{ik}^2]\mathbb{E}[\bar{A}_{jk}^2]\\
&=&\frac{1}{n^2}\sum_{i\neq j\neq k}(a_i^2+a_ia_j+a_ia_k)\mu_{ij}(1-\mu_{ij})\mu_{ik}(1-\mu_{ik})\mu_{jk}(1-\mu_{jk})\\
&=&\Theta\left(\frac{1}{n^2(np_n)}\right).
\end{eqnarray*}
Similarly, the order of the second moment of the second term of (\ref{neq10}) is equal to 
\begin{eqnarray}\label{neq11}
\mathbb{E}\left[\left(\frac{1}{n}\sum_{i\neq j\neq k}\mathbb{E}\left[\frac{1}{d_i(d_i-1)}\right]\bar{A}_{ij}\bar{A}_{ik}\mu_{jk}\right)^2\right]=\Theta\left(\frac{p_n}{n^2(np_n)}\right).
\end{eqnarray}
The order of the second moment of the 3rd term and the 4th term in (\ref{neq10}) are the same as (\ref{neq11}).

The order of the second moment of the 5th term of (\ref{neq10}) is
\begin{eqnarray}\label{neq12}
\mathbb{E}\left[\left(\frac{1}{n}\sum_{i\neq j\neq k}\mathbb{E}\left[\frac{1}{d_i(d_i-1)}\right]\bar{A}_{jk}\mu_{ij}\mu_{ik}\right)^2\right]=\Theta\left(\frac{p_n}{n^2}\right).
\end{eqnarray}
The order of the second moment of the last two terms in (\ref{neq10}) are the same as (\ref{neq12}).

Recall that $p_n=n^{-\alpha}$. Note that
\[\frac{1}{n^2(np_n)}\frac{n^2}{p_n}=\frac{1}{n^{1-2\alpha}}.\]
For $\alpha<\frac{1}{2}$, $\frac{1}{n^2(np_n)}=o\left(\frac{p_n}{n^2}\right)$. For $\alpha>\frac{1}{2}$, $\frac{p_n}{n^2}=o\left(\frac{1}{n^2(np_n)}\right)$. For $\alpha=\frac{1}{2}$, $\frac{1}{n^2(np_n)}=\frac{p_n}{n^2}$. The leading terms of $Y_{1n}$ can be expressed as follows.

If $\alpha>\frac{1}{2}$, then
\begin{eqnarray}\label{y1ng}
Y_{1n}
&=&\frac{1}{n}\sum_{i\neq j\neq k}a_i\bar{A}_{ij}\bar{A}_{ik}\bar{A}_{jk}+o_P\left(\frac{1}{n\sqrt{np_n}}\right).
\end{eqnarray}

If $\alpha<\frac{1}{2}$, then
\begin{eqnarray} \label{y1ns}
Y_{1n}
=\frac{1}{n}\sum_{i\neq j\neq k}a_i\left(\bar{A}_{jk}\mu_{ij}\mu_{ik}+\bar{A}_{ij}\mu_{jk}\mu_{ik}+\bar{A}_{ik}\mu_{jk}\mu_{ij}\right)+o_P\left(\frac{\sqrt{p_n}}{n}\right).
\end{eqnarray}

If $\alpha=\frac{1}{2}$, then
\begin{eqnarray} \nonumber
Y_{1n}
&=&\frac{1}{n}\sum_{i\neq j\neq k}a_i\left(\bar{A}_{ij}\bar{A}_{ik}\bar{A}_{jk}+\bar{A}_{jk}\mu_{ij}\mu_{ik}+\bar{A}_{ij}\mu_{jk}\mu_{ik}+\bar{A}_{ik}\mu_{jk}\mu_{ij}\right)+o_P\left(\frac{\sqrt{p_n}}{n}\right).
\end{eqnarray}

\medskip

{\bf (c) Order of $Y_{2n}$.} Now we show that
\begin{eqnarray}\label{y2n}
Y_{2n}=o_P\left(\frac{\sqrt{p_n}}{n}+\frac{1}{n\sqrt{np_n}}\right).
\end{eqnarray}

Note that
\begin{eqnarray*}
&&\mathbb{E}\left[\big|h^{(k_0)}(X_i)(t_i-\mathbb{E}[t_i])(d_i-\mu_i)^{k_0}\big|\right]\\
&=&\mathbb{E}\left[\big|h^{(k_0)}(X_i)(t_i-\mathbb{E}[t_i])(d_i-\mu_i)^{k_0}\big|I[X_i\leq \delta_nnp_n]\right]\\
&&+\mathbb{E}\left[\big|h^{(k_0)}(X_i)(t_i-\mathbb{E}[t_i])(d_i-\mu_i)^{k_0}\big|I[X_i> \delta_nnp_n]\right].
\end{eqnarray*}
By (\ref{neq15}), it is easy to verify that
\[\mathbb{E}[(t_i-\mathbb{E}[t_i])^2]=O\left(n^3p_n^5+n^2p_n^3\right).\]
By (\ref{neq7}), we get
\begin{eqnarray*}
&&\mathbb{E}\left[\big|h^{(k_0)}(X_i)(t_i-\mathbb{E}[t_i])(d_i-\mu_i)^{k_0}\big|I[X_i> \delta_nnp_n]\right]\\
&=&O\left(\frac{\sqrt{\mathbb{E}\left[(d_i-\mu_i)^{2k_0}\right]\mathbb{E}[(t_i-\mathbb{E}[t_i])^2]}}{\big(\delta_nnp_n\big)^{k_0+2}}\right)\\
&=&O\left(\frac{p_n}{\delta_n^{k_0+2}\big(np_n\big)^{\frac{k_0+1}{2}}}+\frac{1}{\sqrt{n}\delta_n^{k_0+2}\big(np_n\big)^{\frac{k_0+1}{2}}}\right).
\end{eqnarray*}
Since $k_0=\lceil 1+\frac{1}{1-\alpha}\rceil +2> \max\{3,\frac{1}{1-\alpha}\}$, then
\[\frac{p_nn\sqrt{np_n}}{\big(np_n\big)^{\frac{k_0+1}{2}}}=\frac{1}{\big(np_n\big)^{\frac{k_0}{2}-1}}=o(1).\]
\[\frac{n\sqrt{np_n}}{\sqrt{n}\big(np_n\big)^{\frac{k_0+1}{2}}}=o(1).\]
Hence
\begin{eqnarray*}
\mathbb{E}\left[\big|h^{(k_0)}(X_i)(t_i-\mathbb{E}[t_i])(d_i-\mu_i)^{k_0}\big|I[X_i> \delta_nnp_n]\right]
=o\left(\frac{\sqrt{p_n}}{n}+\frac{1}{n\sqrt{np_n}}\right).
\end{eqnarray*}

If $X_i\leq \delta_nnp_n$, then $d_i\leq \delta_nnp_n$. Otherwise, $X_i$ can not be between $d_i$ and $\mu_i$. In this case, $|h^{(k)}(x)|=O(1)$. Then  by Lemma \ref{lem1}, we have
\begin{eqnarray*}
&&\mathbb{E}\left[\big|h^{(k_0)}(X_i)(t_i-\mathbb{E}[t_i])(d_i-\mu_i)^{k_0}\big|I[X_i\leq \delta_nnp_n]\right]\\
&\leq&O(1)\mathbb{E}\left[\big|(t_i-\mathbb{E}[t_i])(d_i-\mu_i)^{k_0}\big|I[d_i\leq \delta_nnp_n]\right]\\
&\leq&O(n^2(np_n)^{k_0}))\sum_{t=2}^{\delta_nnp_n}\mathbb{P}(d_i=t)\\
&=&e^{-np_n\beta(1+o(1))}.
\end{eqnarray*}

Hence,
\begin{eqnarray*}
&&\left|\frac{1}{n}\sum_{i=1}^n(t_i-\mathbb{E}[t_i])\left(\frac{h^{(k_0)}(X_i)}{k_0!}(d_i-\mu_i)^{k_0}-\mathbb{E}\left[\frac{h^{(k_0)}(X_i)}{k_0!}(d_i-\mu_i)^{k_0}\right]\right)\right|\\
&=&o_P\left(\frac{\sqrt{p_n}}{n}+\frac{1}{n\sqrt{np_n}}\right).
\end{eqnarray*}

By (\ref{neq15}) and the definition of $Y_{2n}$, we only need to show the following equations.
\begin{eqnarray}\label{neq1}
\frac{1}{n}\sum_{i\neq j\neq k}\frac{2\mu_i-1}{\mu_i^2(\mu_i-1)^2}(d_i-\mu_i)\bar{A}_{ij}\bar{A}_{ik}\bar{A}_{jk}&=&o_P\left(\frac{\sqrt{p_n}}{n}+\frac{1}{n\sqrt{np_n}}\right) ,  \\ \label{neq2}
\sum_{s=2}^{k_0-1}\frac{1}{n}\sum_{i\neq j\neq k}\frac{h^{(s)}(\mu_i)}{s!}\bar{A}_{ij}\bar{A}_{ik}\bar{A}_{jk}\Big[(d_i-\mu_i)^s-\mathbb{E}[(d_i-\mu_i)^s]\Big] 
&=&o_P\left(\frac{\sqrt{p_n}}{n}+\frac{1}{n\sqrt{np_n}}\right),
\end{eqnarray}
\begin{eqnarray}\label{neq3}
\frac{1}{n}\sum_{i\neq j\neq k}\frac{2\mu_i-1}{\mu_i^2(\mu_i-1)^2}(d_i-\mu_i)\bar{A}_{ij}\bar{A}_{ik}\mu_{jk}&=&o_P\left(\frac{\sqrt{p_n}}{n}+\frac{1}{n\sqrt{np_n}}\right),  \\ \label{neq17}
\sum_{s=2}^{k_0-1}\frac{1}{n}\sum_{i\neq j\neq k}\frac{h^{(s)}(\mu_i)}{s!}\bar{A}_{ij}\bar{A}_{ik}\mu_{jk}\Big[(d_i-\mu_i)^s-\mathbb{E}[(d_i-\mu_i)^s]\Big] 
&=&o_P\left(\frac{\sqrt{p_n}}{n}+\frac{1}{n\sqrt{np_n}}\right),
\end{eqnarray}
\begin{eqnarray}\label{sunneq3}
\frac{1}{n}\sum_{i\neq j\neq k}\frac{2\mu_i-1}{\mu_i^2(\mu_i-1)^2}(d_i-\mu_i)\bar{A}_{ij}\bar{A}_{jk}\mu_{ik}&=&o_P\left(\frac{\sqrt{p_n}}{n}+\frac{1}{n\sqrt{np_n}}\right),  \\ \label{suneq17}
\sum_{s=2}^{k_0-1}\frac{1}{n}\sum_{i\neq j\neq k}\frac{h^{(s)}(\mu_i)}{s!}\bar{A}_{ij}\bar{A}_{jk}\mu_{ik}\Big[(d_i-\mu_i)^s-\mathbb{E}[(d_i-\mu_i)^s]\Big] 
&=&o_P\left(\frac{\sqrt{p_n}}{n}+\frac{1}{n\sqrt{np_n}}\right),
\end{eqnarray}
\begin{eqnarray}\nonumber
&&\frac{1}{n}\sum_{i\neq j\neq k}\frac{2\mu_i-1}{\mu_i^2(\mu_i-1)^2}(d_i-\mu_i)\left(\bar{A}_{ij}\mu_{ik}\mu_{jk}+\bar{A}_{ik}\mu_{ij}\mu_{jk}+\bar{A}_{jk}\mu_{ij}\mu_{ik}\right)\\ \label{neq18}
&&-\frac{2}{n}\sum_{i\neq j\neq k}\frac{2\mu_i-1}{\mu_i^2(\mu_i-1)^2}\mathbb{E}[\bar{A}_{ij}^2]\mu_{ik}\mu_{jk}=o_P\left(\frac{\sqrt{p_n}}{n}+\frac{1}{n\sqrt{np_n}}\right), \\ \nonumber
&&\sum_{s=2}^{k_0-1}\frac{1}{n}\sum_{i\neq j\neq k}\frac{h^{(s)}(\mu_i)}{s!}\left(\bar{A}_{ij}\mu_{ik}\mu_{jk}+\bar{A}_{ik}\mu_{ij}\mu_{jk}+\bar{A}_{jk}\mu_{ij}\mu_{ik}\right)\Big[(d_i-\mu_i)^s-\mathbb{E}[(d_i-\mu_i)^s]\Big]  \\ \label{neq19}
&&=o_P\left(\frac{\sqrt{p_n}}{n}+\frac{1}{n\sqrt{np_n}}\right).
\end{eqnarray}

 The left-hand side of  (\ref{neq1}) can be expressed as
\begin{eqnarray}\nonumber
\frac{1}{n}\sum_{i\neq j\neq k}\frac{2\mu_i-1}{\mu_i^2(\mu_i-1)^2}(d_i-\mu_i)\bar{A}_{ij}\bar{A}_{ik}\bar{A}_{jk} 
&=&\frac{1}{n}\sum_{i\neq j\neq k, l\neq i}\frac{2\mu_i-1}{\mu_i^2(\mu_i-1)^2}\bar{A}_{il}\bar{A}_{ij}\bar{A}_{ik}\bar{A}_{jk}\\ \nonumber
&=&\frac{1}{n}\sum_{i\neq j\neq k\neq l}\frac{2\mu_i-1}{\mu_i^2(\mu_i-1)^2}\bar{A}_{il}\bar{A}_{ij}\bar{A}_{ik}\bar{A}_{jk}\\ \label{neq4}
&&+\frac{2}{n}\sum_{i\neq j\neq k}\frac{2\mu_i-1}{\mu_i^2(\mu_i-1)^2}\bar{A}_{ij}^2\bar{A}_{ik}\bar{A}_{jk}.
\end{eqnarray}
The second moment of the first term in (\ref{neq4}) is bounded by
\begin{eqnarray*}
\mathbb{E}\left[\left(\frac{1}{n}\sum_{i\neq j\neq k\neq l}\frac{2\mu_i-1}{\mu_i^2(\mu_i-1)^2}\bar{A}_{il}\bar{A}_{ij}\bar{A}_{ik}\bar{A}_{jk}\right)^2\right]
=O\left(\frac{n^4p_n^4}{n^2(np_n)^6}\right)
=O\left(\frac{1}{n^2(np_n)^2}\right).
\end{eqnarray*}
The second moment of the second term in (\ref{neq4})  is bounded by
\begin{eqnarray*}
\mathbb{E}\left[\left(\frac{1}{n}\sum_{i\neq j\neq k}\frac{2\mu_i-1}{\mu_i^2(\mu_i-1)^2}\bar{A}_{ij}^2\bar{A}_{ik}\bar{A}_{jk}\right)^2\right]
=O\left(\frac{n^3p_n^3}{n^2(np_n)^6}\right)
=O\left(\frac{1}{n^2(np_n)^3}\right).
\end{eqnarray*}
 Hence (\ref{neq1})  holds.

Now we prove (\ref{neq2}). It is easy to verify that
\begin{eqnarray*}
\mathbb{E}\left[\left(\frac{1}{n}\sum_{i\neq j\neq k}\frac{h^{(s)}(\mu_i)}{s!}\bar{A}_{ij}\bar{A}_{ik}\bar{A}_{jk}\mathbb{E}[(d_i-\mu_i)^s]\right)^2\right]=O\left(\frac{1}{n^2(np_n)^{s+1}}\right).
\end{eqnarray*}

Recall the expression of $(d_i-\mu_i)^s$  in (\ref{dseq}).
Given $\lambda_{t1},\lambda_{t2},\dots,\lambda_{tt}$  such that $\lambda_{t1}+\lambda_{t2}+\dots+\lambda_{tt}=s$. 
Suppose $\lambda_l\geq2$ for all $l=1,2,\dots,t$. If $s\geq3$, then $s-t\geq2$. In this case,
\begin{eqnarray}\nonumber
\mathbb{E}\left[\left|\frac{1}{n}\sum_{j_1\neq j_2\neq \dots,j_t\neq i\neq j\neq k}h^{(s)}(\mu_i)\bar{A}_{ij}\bar{A}_{ik}\bar{A}_{jk}\prod_{l=1}^t\bar{A}_{ij_l}^{\lambda_{l}}\right|\right]
&=&O\left(\frac{(np_n)^{t+3}}{n(np_n)^{s+2}}\right)\\ \nonumber
&=&O\left(\frac{1}{n(np_n)^{s-t-1}}\right)\\ \label{suneq1}
&=&o\left(\frac{1}{n\sqrt{np_n}}+\frac{\sqrt{p_n}}{n}\right).
\end{eqnarray}
When $j\in\{j_1,j_2,\dots,j_t\}$ or $k\in\{j_1,j_2,\dots,j_t\}$ in the summation of (\ref{suneq1}), (\ref{suneq1}) still holds.

If $s=2$, then $t=1$. In this case,
\begin{eqnarray*}
\mathbb{E}\left[\left(\frac{1}{n}\sum_{j_1\neq  i\neq j\neq k}h^{(s)}(\mu_i)\bar{A}_{ij}\bar{A}_{ik}\bar{A}_{jk}\bar{A}_{ij_1}^2\right)^2\right]
&=&O\left(\frac{(np_n)^{5}}{n^2(np_n)^{8}}\right)\\
&=&o\left(\frac{1}{n^2(np_n)}+\frac{p_n}{n^2}\right),\\
\mathbb{E}\left[\left|\frac{1}{n}\sum_{i\neq j\neq k}h^{(s)}(\mu_i)\bar{A}_{ij}^3\bar{A}_{ik}\bar{A}_{jk}\right|\right]
&=&O\left(\frac{(np_n)^{3}}{n(np_n)^{4}}\right)\\
&=&o\left(\frac{1}{n\sqrt{np_n}}+\frac{\sqrt{p_n}}{n}\right).
\end{eqnarray*}
Suppose $\lambda_{t1}=\dots=\lambda_{tt_1}=1$ and $\lambda_{t,t_1+1}\geq2$, \dots, $\lambda_{tt}\geq2$ for $t_1\geq1$. If $s-t\geq2$, then (\ref{suneq1}) holds. Assume $s=t$. In this case, we have
\begin{eqnarray*}
\mathbb{E}\left[\left(\frac{1}{n}\sum_{j_1\neq\dots\neq j_s\neq i\neq j\neq k}h^{(s)}(\mu_i)\bar{A}_{ij}\bar{A}_{ik}\bar{A}_{jk}\bar{A}_{ij_1}\dots\bar{A}_{ij_s}\right)^2\right]
&=&O\left(\frac{(np_n)^{s+3}}{n^2(np_n)^{2s+4}}\right)\\
&=&o\left(\frac{p_n}{n^2}+\frac{p_n}{n^2(np_n)}\right),\\
\mathbb{E}\left[\left(\frac{1}{n}\sum_{j_1\neq\dots\neq j_s\neq i\neq k}h^{(s)}(\mu_i)\bar{A}_{ij_1}^2\bar{A}_{ik}\bar{A}_{j_1k}\bar{A}_{ij_2}\dots\bar{A}_{ij_s}\right)^2\right]
&=&O\left(\frac{(np_n)^{s+2}}{n^2(np_n)^{2s+4}}\right)\\
&=&o\left(\frac{p_n}{n^2}+\frac{p_n}{n^2(np_n)}\right),\\
\mathbb{E}\left[\left|\frac{1}{n}\sum_{j_1\neq\dots\neq j_s\neq i}h^{(s)}(\mu_i)\bar{A}_{ij_1}^2\bar{A}_{ij_2}^2\bar{A}_{j_1j_2}\bar{A}_{ij_3}\dots\bar{A}_{ij_s}\right|\right]
&=&O\left(\frac{(np_n)^{s+1}}{n(np_n)^{s+2}}\right)\\
&=&o\left(\frac{1}{n\sqrt{np_n}}+\frac{\sqrt{p_n}}{n}\right).
\end{eqnarray*}
Assume $t=s-1$. In this case, we have
\begin{eqnarray*}
\mathbb{E}\left[\left(\frac{1}{n}\sum_{j_1\neq\dots\neq j_{s-1}\neq i\neq j\neq k}h^{(s)}(\mu_i)\bar{A}_{ij}\bar{A}_{ik}\bar{A}_{jk}\bar{A}_{ij_1}\dots\bar{A}_{ij_{s-2}}\bar{A}_{ij_{s-1}}^2\right)^2\right]
&=&O\left(\frac{(np_n)^{s+3}}{n^2(np_n)^{2s+4}}\right)\\
&=&o\left(\frac{p_n}{n^2}+\frac{p_n}{n^2(np_n)}\right),\\
\mathbb{E}\left[\left|\frac{1}{n}\sum_{j_1\neq\dots\neq j_{s-1}\neq i\neq k}h^{(s)}(\mu_i)\bar{A}_{ij_1}^2\bar{A}_{ik}\bar{A}_{j_1k}\bar{A}_{ij_2}\dots\bar{A}_{ij_{s-2}}\bar{A}_{ij_{s-1}}^2\right|\right]
&=&O\left(\frac{(np_n)^{s+1}}{n(np_n)^{s+2}}\right)\\
&=&o\left(\frac{1}{n\sqrt{np_n}}+\frac{\sqrt{p_n}}{n}\right),\\
\mathbb{E}\left[\left|\frac{1}{n}\sum_{j_1\neq\dots\neq j_{s-1}\neq i}h^{(s)}(\mu_i)\bar{A}_{ij_1}^2\bar{A}_{ij_2}^2\bar{A}_{j_1j_2}\bar{A}_{ij_3}\dots\bar{A}_{ij_{s-2}}\bar{A}_{ij_{s-1}}^2\right|\right]
&=&O\left(\frac{(np_n)^{s}}{n(np_n)^{s+2}}\right)\\
&=&o\left(\frac{1}{n\sqrt{np_n}}+\frac{\sqrt{p_n}}{n}\right).
\end{eqnarray*}
Hence (\ref{neq2}) holds.

Now we prove (\ref{neq3}).
\begin{eqnarray*} 
\frac{1}{n}\sum_{i\neq j\neq k}\frac{2\mu_i-1}{\mu_i^2(\mu_i-1)^2}(d_i-\mu_i)\bar{A}_{ij}\bar{A}_{ik}\mu_{jk}
&=&\frac{1}{n}\sum_{i\neq j\neq k\neq l}\frac{2\mu_i-1}{\mu_i^2(\mu_i-1)^2}\bar{A}_{il}\bar{A}_{ij}\bar{A}_{ik}\mu_{jk}\\
&&+\frac{2}{n}\sum_{i\neq j\neq k}\frac{2\mu_i-1}{\mu_i^2(\mu_i-1)^2}\bar{A}_{ij}^2\bar{A}_{ik}\mu_{jk},
\end{eqnarray*}
It is easy to verify that
\begin{eqnarray*} 
&&\mathbb{E}\left[\left(\frac{1}{n}\sum_{i\neq j\neq k\neq l}\frac{2\mu_i-1}{\mu_i^2(\mu_i-1)^2}\bar{A}_{il}\bar{A}_{ij}\bar{A}_{ik}\mu_{jk}\right)^2\right]=O\left(\frac{p_n}{n^2(np_n)^2}\right),\\
&&\mathbb{E}\left[\left(\frac{1}{n}\sum_{i\neq j\neq k}\frac{2\mu_i-1}{\mu_i^2(\mu_i-1)^2}\bar{A}_{ij}^2\bar{A}_{ik}\mu_{jk}\right)^2\right]=O\left(\frac{p_n}{n^2(np_n)^2}\right).
\end{eqnarray*}
Hence (\ref{neq3}) holds. Similarly, (\ref{sunneq3}) holds.

Now we prove (\ref{neq17}). 
\begin{eqnarray*}
\mathbb{E}\left[\left(\frac{1}{n}\sum_{i\neq j\neq k}h^{(s)}(\mu_i)\bar{A}_{ij}\bar{A}_{ik}\mu_{jk}\mathbb{E}[(d_i-\mu_i)^s]\right)^2\right]
&=&O\left(\frac{p_n}{n^2(np_n)^{s+1}}\right)\\
&=&o\left(\frac{1}{n^2np_n}+\frac{p_n}{n^2}\right).
\end{eqnarray*}

Suppose $\lambda_{tl}\geq2$ for all $l=1,2,\dots,t$. If $s-t\geq2$, that is, $s\geq3$, then
\begin{eqnarray}\nonumber
\mathbb{E}\left[\left|\frac{1}{n}\sum_{j_1\neq j_2\neq \dots\neq j_t\neq i\neq j\neq k}h^{(s)}(\mu_i)\bar{A}_{ij}\bar{A}_{ik}\mu_{jk}\prod_{l=1}^t\bar{A}_{ij_l}^{\lambda_{tl}}\right|\right]
&=&O\left(\frac{n^{t+3}p_n^{t+3}}{n(np_n)^{s+2}}\right)\\ \nonumber
&=&O\left(\frac{1}{n(np_n)^{s-t-1}}\right)\\ \label{suneq3}
&=&o\left(\frac{1}{n\sqrt{np_n}}+\frac{\sqrt{p_n}}{n}\right)
\end{eqnarray}
If $s=2$, then $t=1$ and
\begin{eqnarray*}
\mathbb{E}\left[\left(\frac{1}{n}\sum_{j_1\neq i\neq j\neq k}h^{(2)}(\mu_i)\bar{A}_{ij}\bar{A}_{ik}\mu_{jk}\bar{A}_{ij_1}^{2}\right)^2\right]
&=&O\left(\frac{n^{5}p_n^{6}}{n^2(np_n)^{8}}\right)=o\left(\frac{1}{n^2np_n}+\frac{p_n}{n^2}\right),
\end{eqnarray*}
and
\begin{eqnarray*}
\mathbb{E}\left[\left|\frac{1}{n}\sum_{i\neq j\neq k}h^{(2)}(\mu_i)\bar{A}_{ik}\mu_{jk}\bar{A}_{ij}^{3}\right|\right]
&=&O\left(\frac{n^{3}p_n^{3}}{n(np_n)^{4}}\right)=o\left(\frac{1}{n\sqrt{np_n}}+\frac{\sqrt{p_n}}{n}\right).
\end{eqnarray*}

Suppose $\lambda_{t1}=\dots=\lambda_{tt_1}=1$ and $\lambda_{t,t_1+1}\geq2$, $\lambda_{tt}\geq2$ for $t_1\geq1$. If $s-t\geq2$, then (\ref{suneq3}) holds. Assume $s=t$. Then
\begin{eqnarray}\nonumber
\mathbb{E}\left[\left(\frac{1}{n}\sum_{j_1\neq\dots\neq j_s\neq i\neq j\neq k}h^{(s)}(\mu_i)\bar{A}_{ij}\bar{A}_{ik}\mu_{jk}\bar{A}_{ij_1}\dots\bar{A}_{ij_s}\right)^2\right]
&=&O\left(\frac{n^{s+3}p_n^{s+4}}{n^2(np_n)^{2s+4}}\right)\\ \label{suneq10}
&=&o\left(\frac{1}{n^2np_n}+\frac{p_n}{n^2}\right).
\end{eqnarray}
For $j\in\{j_1,j_2,\dots,j_s\}$ or $k\in\{j_1,j_2,\dots,j_s\}$ in the summation of (\ref{suneq10}), (\ref{suneq10}) still hods.

Assume $t=s-1$. Then
\begin{eqnarray}\nonumber
\mathbb{E}\left[\left(\frac{1}{n}\sum_{j_1\neq\dots\neq j_{s-1}\neq i\neq j\neq k}h^{(s)}(\mu_i)\bar{A}_{ij}\bar{A}_{ik}\mu_{jk}\bar{A}_{ij_1}\dots\bar{A}_{ij_{s-2}}\bar{A}_{ij_{s-1}}^2\right)^2\right]
&=&O\left(\frac{n^{s+3}p_n^{s+3}}{n^2(np_n)^{2s+4}}\right)\\ \nonumber\label{suneq11}
&=&o\left(\frac{1}{n^2np_n}+\frac{p_n}{n^2}\right).
\end{eqnarray}
\begin{eqnarray}\nonumber
\mathbb{E}\left[\left|\frac{1}{n}\sum_{\substack{j_1\neq\dots\neq j_{s-1}\neq i\\ \{j,k\}\cap\{j_1,\dots,j_{s-1}\}\neq\emptyset}}h^{(s)}(\mu_i)\bar{A}_{ij}\bar{A}_{ik}\mu_{jk}\bar{A}_{ij_1}\dots\bar{A}_{ij_{s-2}}\bar{A}_{ij_{s-1}}^2\right|\right]
&=&O\left(\frac{n^{s+1}p_n^{s+1}}{n(np_n)^{s+2}}\right)\\ \label{suneq11}
&=&o\left(\frac{1}{n^2np_n}+\frac{p_n}{n^2}\right).
\end{eqnarray}
Hence (\ref{neq17}) holds. Similarly (\ref{suneq17}) holds.

Now we prove (\ref{neq18}).
\begin{eqnarray}\nonumber
&&\frac{1}{n}\sum_{i\neq j\neq k}\frac{2\mu_i-1}{\mu_i^2(\mu_i-1)^2}(d_i-\mu_i)\left(\bar{A}_{ij}\mu_{ik}\mu_{jk}+\bar{A}_{ik}\mu_{ij}\mu_{jk}\right)-\frac{2}{n}\sum_{i\neq j\neq k}\frac{2\mu_i-1}{\mu_i^2(\mu_i-1)^2}\mathbb{E}[\bar{A}_{ij}^2]\mu_{ik}\mu_{jk} \\ \nonumber
&=&\frac{1}{n}\sum_{i\neq j\neq k\neq l}\frac{2\mu_i-1}{\mu_i^2(\mu_i-1)^2}\bar{A}_{il}\bar{A}_{ij}\mu_{ik}\mu_{jk}+\frac{1}{n}\sum_{i\neq j\neq k\neq l}\frac{2\mu_i-1}{\mu_i^2(\mu_i-1)^2}\bar{A}_{il}\bar{A}_{ik}\mu_{ij}\mu_{jk}\\
&&+\frac{2}{n}\sum_{i\neq j\neq k}\frac{2\mu_i-1}{\mu_i^2(\mu_i-1)^2}\left(\bar{A}_{ij}^2-\mathbb{E}[\bar{A}_{ij}^2]\right)\mu_{ik}\mu_{jk}\\
&&+\frac{1}{n}\sum_{i\neq j\neq k}\frac{2\mu_i-1}{\mu_i^2(\mu_i-1)^2}\bar{A}_{ik}\bar{A}_{ij}\mu_{ik}\mu_{jk}+\frac{1}{n}\sum_{i\neq j\neq k}\frac{2\mu_i-1}{\mu_i^2(\mu_i-1)^2}\bar{A}_{ij}\bar{A}_{ik}\mu_{ij}\mu_{jk}
\end{eqnarray}
It is easy to verify that 
\begin{eqnarray*}\nonumber
&&\mathbb{E}\left[\left(\frac{1}{n}\sum_{i\neq j\neq k\neq l}\frac{2\mu_i-1}{\mu_i^2(\mu_i-1)^2}\bar{A}_{il}\bar{A}_{ij}\mu_{ik}\mu_{jk}\right)^2\right]=O\left(\frac{p_n}{n^2(np_n)}\right),\\
&&\mathbb{E}\left[\left(\frac{1}{n}\sum_{i\neq j\neq k}\frac{2\mu_i-1}{\mu_i^2(\mu_i-1)^2}\left(\bar{A}_{ij}^2-\mathbb{E}[\bar{A}_{ij}^2]\right)\mu_{ik}\mu_{jk}\right)^2\right]=O\left(\frac{p_n}{n^2(np_n)^2}\right),\\
&&\mathbb{E}\left[\left(\frac{1}{n}\sum_{i\neq j\neq k}\frac{2\mu_i-1}{\mu_i^2(\mu_i-1)^2}\bar{A}_{ik}\bar{A}_{ij}\mu_{ik}\mu_{jk}\right)^2\right]=O\left(\frac{p_n^6}{n^2(np_n)^3}\right).
\end{eqnarray*}

Moreover, we have
\begin{eqnarray}\nonumber
\frac{1}{n}\sum_{i\neq j\neq k}\frac{2\mu_i-1}{\mu_i^2(\mu_i-1)^2}(d_i-\mu_i)\bar{A}_{jk}\mu_{ij}\mu_{ik}&=&\frac{1}{n}\sum_{i\neq j\neq k\neq s}\frac{2\mu_i-1}{\mu_i^2(\mu_i-1)^2}\bar{A}_{is}\bar{A}_{jk}\mu_{ij}\mu_{ik}\\ \nonumber
&&+\frac{2}{n}\sum_{i\neq j\neq k}\frac{2\mu_i-1}{\mu_i^2(\mu_i-1)^2}\bar{A}_{ij}\bar{A}_{jk}\mu_{ij}\mu_{ik},
\end{eqnarray}
and
\begin{eqnarray}\nonumber
\mathbb{E}\left[\left(\frac{1}{n}\sum_{i\neq j\neq k\neq s}\frac{2\mu_i-1}{\mu_i^2(\mu_i-1)^2}\bar{A}_{is}\bar{A}_{jk}\mu_{ij}\mu_{ik}\right)^2\right]&=&O\left(\frac{p_n^2}{n^2(np_n)^2}\right),\\ \nonumber
\mathbb{E}\left[\left(\frac{1}{n}\sum_{i\neq j\neq k}\frac{2\mu_i-1}{\mu_i^2(\mu_i-1)^2}\bar{A}_{ij}\bar{A}_{jk}\mu_{ij}\mu_{ik}\right)^2\right]&=&O\left(\frac{p_n^3}{n^2(np_n)^3}\right).
\end{eqnarray}
Then (\ref{neq18}) holds.

Now we prove (\ref{neq19}). 
\begin{eqnarray}\nonumber
\mathbb{E}\left[\left(\frac{1}{n}\sum_{i\neq j\neq k}\frac{h^{(s)}(\mu_i)}{s!}\bar{A}_{ij}\mu_{ik}\mu_{jk}\mathbb{E}[(d_i-\mu_i)^s] \right)^2\right] 
=O\left(\frac{p_n}{n^2(np_n)^s}\right).
\end{eqnarray}

Suppose $\lambda_l\geq2$ for all $l=1,2,\dots,t$. If $s\geq3$, then
\begin{eqnarray}\nonumber
\mathbb{E}\left[\left|\frac{1}{n}\sum_{j_1\neq j_2\neq \dots,j_t\neq i\neq j\neq k}h^{(s)}(\mu_i)\bar{A}_{ij}\mu_{ik}\mu_{jk}\prod_{l=1}^t\bar{A}_{ij_l}^{\lambda_{l}}\right|\right]
&=&O\left(\frac{n^{t+3}p_n^{t+3}}{n(np_n)^{s+2}}\right)\\ \nonumber
&=&O\left(\frac{1}{n(np_n)^{s-t-1}}\right)\\ \label{suneq4}
&=&o\left(\frac{1}{n\sqrt{np_n}}+\frac{\sqrt{p_n}}{n}\right).
\end{eqnarray}
If $s=2$, then $t=1$ and
\begin{eqnarray*}
\mathbb{E}\left[\left(\frac{1}{n}\sum_{j_1\neq  i\neq j\neq k}h^{(s)}(\mu_i)\bar{A}_{ij}\mu_{ik}\mu_{jk}\bar{A}_{ij_1}^2\right)^2\right]
&=&O\left(\frac{n^{6}p_n^{7}}{n^2(np_n)^{8}}\right)\\
&=&o\left(\frac{1}{n\sqrt{np_n}}+\frac{\sqrt{p_n}}{n}\right),
\end{eqnarray*}

\begin{eqnarray*}
\mathbb{E}\left[\left|\frac{1}{n}\sum_{ i\neq j\neq k}h^{(s)}(\mu_i)\bar{A}_{ij}\mu_{ik}\mu_{jk}\bar{A}_{ik}^2\right|\right]
&=&O\left(\frac{n^{3}p_n^{4}}{n(np_n)^{4}}\right)\\
&=&o\left(\frac{1}{n\sqrt{np_n}}+\frac{\sqrt{p_n}}{n}\right),
\end{eqnarray*}
\begin{eqnarray*}
\mathbb{E}\left[\left|\frac{1}{n}\sum_{i\neq j\neq k}h^{(s)}(\mu_i)\bar{A}_{ij}^3\mu_{ik}\mu_{jk}\right|\right]
&=&O\left(\frac{n^{3}p_n^{3}}{n(np_n)^{4}}\right)\\
&=&o\left(\frac{1}{n\sqrt{np_n}}+\frac{\sqrt{p_n}}{n}\right).
\end{eqnarray*}

Suppose $\lambda_{t1}=\dots=\lambda_{tt_1}=1$ and $\lambda_{t,t_1+1}\geq2$, \dots, $\lambda_{tt}\geq2$ for $t_1\geq1$. If $s-t\geq2$, then (\ref{suneq4}) holds. Assume $s=t$. Then
\begin{eqnarray}\nonumber
\mathbb{E}\left[\left(\frac{1}{n}\sum_{j_1\neq \dots \neq j_s\neq  j\neq i\neq k}h^{(s)}(\mu_i)\bar{A}_{ij}\mu_{ik}\mu_{jk}\bar{A}_{ij_1}\dots\bar{A}_{ij_s}\right)^2\right]&=&O\left(\frac{n^{s+4}p_n^{s+5}}{n^2(np_n)^{2s+4}}\right)\\ \label{seq1}
&=&o\left(\frac{1}{n^2np_n}+\frac{p_n}{n^2}\right),
\end{eqnarray}
and
\begin{eqnarray}\nonumber
\mathbb{E}\left[\left(\frac{1}{n}\sum_{j_1\neq \dots \neq j_s\neq  i\neq k}h^{(s)}(\mu_i)\bar{A}_{ij_1}^2\mu_{ik}\mu_{j_1k}\bar{A}_{ij_2}\dots\bar{A}_{ij_s}\right)^2\right]&=&O\left(\frac{n^{s+4}p_n^{s+4}}{n^2(np_n)^{2s+4}}\right)\\ \label{seq1}
&=&o\left(\frac{1}{n^2np_n}+\frac{p_n}{n^2}\right),
\end{eqnarray}
Assume $t=s-1$. Then
\begin{eqnarray}\nonumber
\mathbb{E}\left[\left(\frac{1}{n}\sum_{j_1\neq \dots \neq j_{s-1}\neq  j\neq i\neq k}h^{(s)}(\mu_i)\bar{A}_{ij}\mu_{ik}\mu_{jk}\bar{A}_{ij_1}\dots\bar{A}_{ij_{s-2}}\bar{A}_{ij_{s-1}}^2\right)^2\right]&=&O\left(\frac{n^{s+4}p_n^{s+5}}{n^2(np_n)^{2s+4}}\right)\\ \label{seq1}
&=&o\left(\frac{1}{n^2np_n}+\frac{p_n}{n^2}\right),
\end{eqnarray}
and
\begin{eqnarray}\nonumber
\mathbb{E}\left[\left|\frac{1}{n}\sum_{\substack{j_1\neq\dots\neq j_{s-1}\neq i\neq k\\ j\in\{j_1,\dots,j_{s-1}\}}}h^{(s)}(\mu_i)\bar{A}_{ij}\mu_{ik}\mu_{jk}\bar{A}_{ij_1}\dots\bar{A}_{ij_{s-2}}\bar{A}_{ij_{s-1}}^2\right|\right]&=&O\left(\frac{n^{s+1}p_n^{s+1}}{n(np_n)^{s+2}}\right)\\ \label{seq1}
&=&o\left(\frac{1}{n^2np_n}+\frac{p_n}{n^2}\right).
\end{eqnarray}
Hence (\ref{neq19}) holds.

By (\ref{eq1}), (\ref{y3n}),(\ref{y1ng}), (\ref{y1ns}), (\ref{y1ne}) and (\ref{y2n}), $\overline{\mathcal{C}}_n-\mathbb{E}[\overline{\mathcal{C}}_n]$ has the following asymptotic expression.
\begin{equation}\label{eqce}
\overline{\mathcal{C}}_n-\mathbb{E}[\overline{\mathcal{C}}_n] =
    \begin{cases}
         \frac{2}{n}\sum_{i< j< k}a_{ijk}\bar{A}_{ij}\bar{A}_{ik}\bar{A}_{jk}+o_P\left(\frac{1}{n\sqrt{np_n}}\right), & \text{if $\alpha>\frac{1}{2}$}, \\
      \frac{2}{n}\sum_{i< j< k}a_{ijk}\bar{A}_{ij}\bar{A}_{ik}\bar{A}_{jk}
+\frac{2}{n}\sum_{i< j}e_{ij}\bar{A}_{ij}+o_P\left(\frac{\sqrt{p_n}}{n}+\frac{1}{n\sqrt{np_n}}\right), & \text{if $\alpha=\frac{1}{2}$}, \\
      \frac{2}{n}\sum_{i< j}e_{ij}\bar{A}_{ij}+o_P\left(\frac{\sqrt{p_n}}{n}\right)  ,& \text{if $\alpha<\frac{1}{2}$}.     
    \end{cases}       
\end{equation}

By Lemma \ref{nlem1}, the proof is complete.

\qed

\subsection{Proof of Theorem \ref{mainthm2}}

For convenience, let $\Delta_{ijk}=A_{ij}A_{jk}A_{ki}$. For distinct indices $i,j,k$, define $d_{i(jk)}=2+\sum_{l\notin \{i,j,k\}}A_{il}$. Then
\begin{eqnarray*}    \mathcal{T}_n=\sum_{i<j<k}\frac{A_{ij}A_{jk}A_{ki}}{d_{i(jk)}d_{j(ki)}d_{k(ij)}}.
\end{eqnarray*}
Note that $d_{i(jk)}$, $d_{j(ki)}$, $d_{k(ij)}$ are independent. With a little abuse of notation, we still write $d_{i(jk)}$ as $d_i$ for convenience. 
Let $g(x,y,z)=\frac{1}{xyz}$. Then the $k$-th derivative of $g$ is 
\[\frac{\partial ^kg(x,y,z)}{\partial x^r\partial y^m\partial z^l}=\frac{r!m!l!(-1)^{r+m+l}}{x^{r+1}y^{m+1}x^{l+1}},\]
where $r,m,l$ are non-negative integers such that $r+m+l=k$. Let $k_0=3+\lceil\frac{1}{1-\alpha}\rceil$. By Taylor expansion, we have
\begin{eqnarray}\nonumber
    \sum_{i<j<k}\frac{\Delta_{ijk}}{d_id_jd_k}&=&\sum_{i<j<k}\frac{\Delta_{ijk}}{\mu_i\mu_j\mu_k}-\sum_{i<j<k}\left(\frac{\Delta_{ijk}(d_i-\mu_i)}{\mu_i^2\mu_j\mu_k}+\frac{\Delta_{ijk}(d_j-\mu_j)}{\mu_i\mu_j^2\mu_k}+\frac{\Delta_{ijk}(d_k-\mu_k)}{\mu_i\mu_j\mu_k^2}\right)\\ \nonumber
    &&+\sum_{t=2}^{k_0-1}\sum_{r+m+l=t}(-1)^t\sum_{i<j<k}\frac{\Delta_{ijk}(d_i-\mu_i)^r(d_j-\mu_j)^m(d_k-\mu_k)^l}{\mu_i^{r+1}\mu_j^{m+1}\mu_k^{l+1}}\\ \label{trieq1}
    &&+(-1)^{k_0}\sum_{r+m+l=k_0}\sum_{i<j<k}\frac{\Delta_{ijk}(d_i-\mu_i)^{r}(d_j-\mu_j)^{m}(d_k-\mu_k)^{l}}{X_i^{r+1}X_j^{m+1}X_k^{l+1}},
\end{eqnarray}
where $X_s$ is between $\mu_s$ and $d_s$ for $s\in\{i,j,k\}$. Then
\begin{eqnarray}\nonumber   &&\sum_{i<j<k}\left(\frac{\Delta_{ijk}}{d_id_jd_k}-\mathbb{E}\left[\frac{\Delta_{ijk}}{d_id_jd_k}\right]\right)\\ \nonumber
    &=&\sum_{i<j<k}\frac{\Delta_{ijk}-\mathbb{E}[\Delta_{ijk}]}{\mu_i\mu_j\mu_k}-\sum_{i<j<k}\left(\frac{\Delta_{ijk}(d_i-\mu_i)}{\mu_i^2\mu_j\mu_k}+\frac{\Delta_{ijk}(d_j-\mu_j)}{\mu_i\mu_j^2\mu_k}+\frac{\Delta_{ijk}(d_k-\mu_k)}{\mu_i\mu_j\mu_k^2}\right)\\ \nonumber
    &&+\sum_{t=2}^{k_0-1}\sum_{r+m+l=t}(-1)^t\sum_{i<j<k}\frac{\Delta_{ijk}(d_i-\mu_i)^r(d_j-\mu_j)^m(d_k-\mu_k)^l}{\mu_i^{r+1}\mu_j^{m+1}\mu_k^{l+1}}\\  \nonumber
    &&-\sum_{t=2}^{k_0-1}\sum_{r+m+l=t}(-1)^t\sum_{i<j<k}\frac{\mathbb{E}[\Delta_{ijk}(d_i-\mu_i)^r(d_j-\mu_j)^m(d_k-\mu_k)^l]}{\mu_i^{r+1}\mu_j^{m+1}\mu_k^{l+1}}\\  \nonumber\label{trieq1}
    &&+(-1)^{k_0}\sum_{r+m+l=k_0}\sum_{i<j<k}\frac{\Delta_{ijk}(d_i-\mu_i)^{r}(d_j-\mu_j)^{m}(d_k-\mu_k)^{l}}{X_i^{r+1}X_j^{m+1}X_k^{l+1}}\\ \label{edeq}
    &&-(-1)^{k_0}\sum_{r+m+l=k_0}\sum_{i<j<k}\mathbb{E}\left[\frac{\Delta_{ijk}(d_i-\mu_i)^{r}(d_j-\mu_j)^{m}(d_k-\mu_k)^{l}}{X_i^{r+1}X_j^{m+1}X_k^{l+1}}\right].
\end{eqnarray}

Next we find the order of each term in (\ref{edeq}).

By (\ref{neq15}), we have
\begin{eqnarray}\nonumber
    \sum_{i<j<k}\frac{\Delta_{ijk}-\mathbb{E}[\Delta_{ijk}]}{\mu_i\mu_j\mu_k}
    &=&\sum_{i<j<k}\frac{\bar{A}_{ij}\bar{A}_{jk}\bar{A}_{ki}}{\mu_i\mu_j\mu_k}+\sum_{i<j<k}\frac{\bar{A}_{ij}\bar{A}_{jk}\mu_{ki}+\bar{A}_{ij}\mu_{jk}\bar{A}_{ki}+\mu_{ij}\bar{A}_{jk}\bar{A}_{ki}}{\mu_i\mu_j\mu_k}\\ \label{trieq2}
&&+\sum_{i<j<k}\frac{\bar{A}_{ij}\mu_{jk}\mu_{ki}+\mu_{ij}\bar{A}_{jk}\mu_{ki}+\mu_{ij}\mu_{jk}\bar{A}_{ki}}{\mu_i\mu_j\mu_k}.
\end{eqnarray}

Now we find the order of each term in (\ref{trieq2}) by calculating the second moment of each term.  Note that
\begin{eqnarray}\nonumber
   \mathbb{E}\left[\left( \sum_{i<j<k}\frac{\bar{A}_{ij}\bar{A}_{jk}\bar{A}_{ki}}{\mu_i\mu_j\mu_k}\right)^2\right]&=&\sum_{i<j<k}\frac{\mathbb{E}[\bar{A}_{ij}^2]\mathbb{E}[\bar{A}_{jk}^2]\mathbb{E}[\bar{A}_{ki}^2]}{\mu_i^2\mu_j^2\mu_k^2}\\ \nonumber
   &=&\sum_{i<j<k}\frac{\mu_{ij}(1-\mu_{ij})\mu_{jk}(1-\mu_{jk})\mu_{ki}(1-\mu_{ki})}{\mu_i^2\mu_j^2\mu_k^2}\\ \label{trieq3}
   &=&\Theta\left(\frac{1}{(np_n)^3}\right),
\end{eqnarray}
\begin{eqnarray}\nonumber
\mathbb{E}\left[\left(\sum_{i<j<k}\frac{\bar{A}_{ij}\bar{A}_{jk}\mu_{ki}}{\mu_i\mu_j\mu_k}\right)^2\right]&=&\sum_{i<j<k}\frac{\mathbb{E}[\bar{A}_{ij}^2]\mathbb{E}[\bar{A}_{jk}^2]\mu_{ki}^2}{\mu_i^2\mu_j^2\mu_k^2}\\ \nonumber
   &=&\sum_{i<j<k}\frac{\mu_{ij}(1-\mu_{ij})\mu_{jk}(1-\mu_{jk})\mu_{ki}^2}{\mu_i^2\mu_j^2\mu_k^2}\\ \label{trieq4}
   &=&\Theta\left(\frac{p_n}{(np_n)^3}\right),
\end{eqnarray}
and
\begin{eqnarray}\nonumber
\mathbb{E}\left[\left(\sum_{i<j<k}\frac{\bar{A}_{ij}\mu_{jk}\mu_{ki}}{\mu_i\mu_j\mu_k}\right)^2\right] &=&\sum_{i<j<k,i<j<k_1}\frac{\mathbb{E}[\bar{A}_{ij}^2]\mu_{jk}\mu_{ki}\mu_{jk_1}\mu_{k_1i}}{\mu_i^2\mu_j^2\mu_k\mu_{k_1}}\\ \nonumber
&=&\sum_{i<j<k,k_1\neq k}\frac{\mu_{ij}(1-\mu_{ij})\mu_{jk}\mu_{ki}\mu_{jk_1}\mu_{k_1i}}{\mu_i^2\mu_j^2\mu_k\mu_{k_1i}}+\sum_{i<j<k}\frac{\mu_{ij}(1-\mu_{ij})\mu_{jk}^2\mu_{ki}^2}{\mu_i^2\mu_j^2\mu_k^2}\\ \label{trieq5}
   &=&\Theta\left(\frac{p_n}{(np_n)^2}\right).
\end{eqnarray}
Combining (\ref{trieq3}), (\ref{trieq4}) and (\ref{trieq5}), we get the leading terms of (\ref{trieq2}) as follows.

If $\alpha>\frac{1}{2}$, then 
\begin{eqnarray}\label{trieq6}
    \sum_{i<j<k}\frac{\Delta_{ijk}-\mathbb{E}[\Delta_{ijk}]}{\mu_i\mu_j\mu_k}
=\sum_{i<j<k}\frac{\bar{A}_{ij}\bar{A}_{jk}\bar{A}_{ki}}{\mu_i\mu_j\mu_k}+o_P\left(\frac{1}{np_n\sqrt{np_n}}\right).
\end{eqnarray}

If $\alpha<\frac{1}{2}$, then 
\begin{eqnarray}\label{trieq7}
    \sum_{i<j<k}\frac{\Delta_{ijk}-\mathbb{E}[\Delta_{ijk}]}{\mu_i\mu_j\mu_k}   &=&\sum_{i<j<k}\frac{\bar{A}_{ij}\mu_{jk}\mu_{ki}+\mu_{ij}\bar{A}_{jk}\mu_{ki}+\mu_{ij}\mu_{jk}\bar{A}_{ki}}{\mu_i\mu_j\mu_k}+o_P\left(\frac{1}{n\sqrt{p_n}}\right).
\end{eqnarray}

If $\alpha=\frac{1}{2}$, then 
\begin{eqnarray}\nonumber
    \sum_{i<j<k}\frac{\Delta_{ijk}-\mathbb{E}[\Delta_{ijk}]}{\mu_i\mu_j\mu_k}   &=&\sum_{i<j<k}\frac{\bar{A}_{ij}\bar{A}_{jk}\bar{A}_{ki}}{\mu_i\mu_j\mu_k}+\sum_{i<j<k}\frac{\bar{A}_{ij}\mu_{jk}\mu_{ki}+\mu_{ij}\bar{A}_{jk}\mu_{ki}+\mu_{ij}\mu_{jk}\bar{A}_{ki}}{\mu_i\mu_j\mu_k}\\ \label{trieq8}
    &&+o_P\left(\frac{1}{n\sqrt{p_n}}+\frac{1}{np_n\sqrt{np_n}}\right).
\end{eqnarray}

Next, we consider the second term of (\ref{edeq}). By (\ref{neq15}), one gets
\begin{eqnarray*}
\sum_{i\neq j\neq k}\frac{\Delta_{ijk}(d_i-\mu_i)}{\mu_i^2\mu_j\mu_k}&=&\sum_{i\neq j\neq k\neq l}\frac{\bar{A}_{ij}\bar{A}_{jk}\bar{A}_{ki}\bar{A}_{il}}{\mu_i^2\mu_j\mu_k}+3\sum_{i\neq j\neq k\neq l}\frac{\bar{A}_{ij}\bar{A}_{jk}\mu_{ki}\bar{A}_{il}}{\mu_i^2\mu_j\mu_k}\\
&&+3\sum_{i\neq j\neq k\neq l}\frac{\bar{A}_{ij}\mu_{jk}\mu_{ki}\bar{A}_{il}}{\mu_i^2\mu_j\mu_k}
+\sum_{i\neq j\neq k\neq l}\frac{\mu_{ij}\mu_{jk}\mu_{ki}\bar{A}_{il}}{\mu_i^2\mu_j\mu_k}.
\end{eqnarray*}
Similar to (\ref{trieq3}), (\ref{trieq4}) and (\ref{trieq5}), it is easy to verify that
\begin{eqnarray*}\nonumber
\sum_{i\neq j\neq k\neq l}\frac{\bar{A}_{ij}\bar{A}_{jk}\bar{A}_{ki}\bar{A}_{il}}{\mu_i^2\mu_j\mu_k}&=&O_P\left(\frac{1}{(np_n)^2}\right),\\ \nonumber
\sum_{i\neq j\neq k\neq l}\frac{\bar{A}_{ij}\bar{A}_{jk}\mu_{ki}\bar{A}_{il}}{\mu_i^2\mu_j\mu_k}&=&O_P\left(\frac{\sqrt{p_n}}{(np_n)^2}\right),\\ \nonumber
\sum_{i\neq j\neq k\neq l}\frac{\bar{A}_{ij}\mu_{jk}\mu_{ki}\bar{A}_{il}}{\mu_i^2\mu_j\mu_k}&=&O_P\left(\frac{1}{np_n\sqrt{n}}\right), \\ 
\sum_{i\neq j\neq k\neq l}\frac{\mu_{ij}\mu_{jk}\mu_{ki}\bar{A}_{il}}{\mu_i^2\mu_j\mu_k}&=&O_P\left(\frac{\sqrt{p_n}}{np_n}\right).
\end{eqnarray*}
Hence the second term of (\ref{edeq}) is equal to
\begin{eqnarray}\nonumber
&&\sum_{i<j<k}\left(\frac{\Delta_{ijk}(d_i-\mu_i)}{\mu_i^2\mu_j\mu_k}+\frac{\Delta_{ijk}(d_j-\mu_j)}{\mu_i\mu_j^2\mu_k}+\frac{\Delta_{ijk}(d_k-\mu_k)}{\mu_i\mu_j\mu_k^2}\right)\\  \nonumber
&=&\sum_{i< j< k, l\notin\{i,j,k\}}\left(\frac{\mu_{ij}\mu_{jk}\mu_{ki}\bar{A}_{il}}{\mu_i^2\mu_j\mu_k}+\frac{\mu_{ij}\mu_{jk}\mu_{ki}\bar{A}_{jl}}{\mu_i\mu_j^2\mu_k}+\frac{\mu_{ij}\mu_{jk}\mu_{ki}\bar{A}_{kl}}{\mu_i\mu_j\mu_k^2}\right)\\ \label{trieq9}
&&+o_P\left(\frac{1}{n\sqrt{p_n}}+\frac{1}{np_n\sqrt{np_n}}\right).
\end{eqnarray}

Now we consider the last two terms of (\ref{edeq}). Given $r,m,l$ such that $r+m+l=k_0$, we have
\begin{eqnarray}\nonumber
&&\sum_{i<j<k}\mathbb{E}\left[\left|\frac{\Delta_{ijk}(d_i-\mu_i)^{r}(d_j-\mu_j)^{m}(d_k-\mu_k)^{l}}{X_i^{r+1}X_j^{m+1}X_k^{l+1}}\right|\right]\\ \nonumber
&=&\sum_{i<j<k}\Bigg(\mathbb{E}\left[\left|\frac{\Delta_{ijk}(d_i-\mu_i)^{r}(d_j-\mu_j)^{m}(d_k-\mu_k)^{l}}{X_i^{r+1}X_j^{m+1}X_k^{l+1}}\right|I[X_i\geq\delta_nnp_n,X_j\geq\delta_nnp_n,X_k\geq\delta_nnp_n]\right]\\  \nonumber
&&+\mathbb{E}\left[\left|\frac{\Delta_{ijk}(d_i-\mu_i)^{r}(d_j-\mu_j)^{m}(d_k-\mu_k)^{l}}{X_i^{r+1}X_j^{m+1}X_k^{l+1}}\right|I[X_i<\delta_nnp_n\ \cup\ X_j<\delta_nnp_n\cup X_k<\delta_nnp_n]\right]\Bigg). \\ \label{trieq10} 
\end{eqnarray}

Recall that $k_0=3+\lceil\frac{1}{1-\alpha}\rceil$. Then $k_0>\max\{3,1+\frac{1}{1-\alpha}\}$. Moreover, $d_i,d_j,d_k$ are independent. Hence,
\begin{eqnarray}\nonumber
&&\sum_{i<j<k}\mathbb{E}\left[\left|\frac{\Delta_{ijk}(d_i-\mu_i)^{r}(d_j-\mu_j)^{m}(d_k-\mu_k)^{l}}{X_i^{r+1}X_j^{m+1}X_k^{l+1}}\right|I[X_i\geq\delta_nnp_n,X_j\geq\delta_nnp_n,X_k\geq\delta_nnp_n]\right]\\ \nonumber
&\leq&\frac{1}{\delta_n^{k_0+3}(np_n)^{k_0+3}}\sum_{i<j<k}\mathbb{E}[\Delta_{ijk}]\mathbb{E}[|(d_i-\mu_i)^{r}|]\mathbb{E}[|(d_j-\mu_j)^{m}|]\mathbb{E}[|(d_k-\mu_k)^{l}|]\\ \nonumber
&=&O\left(\frac{(np_n)^{\frac{k_0}{2}+3}}{\delta_n^{k_0+3}(np_n)^{k_0+3}}\right)\\ \label{trieq11} 
&=&o\left(\frac{1}{np_n\sqrt{np_n}}+\frac{\sqrt{p_n}}{np_n}\right).
\end{eqnarray}
Since $d_i\geq1$ and $X_i$ is between $d_i$ and $\mu_i$, then $X_i\geq1$. Hence
\begin{eqnarray} \nonumber
\sum_{i<j<k}&\mathbb{E}&\left[\left|\frac{\Delta_{ijk}(d_i-\mu_i)^{r}(d_j-\mu_j)^{m}(d_k-\mu_k)^{l}}{X_i^{r+1}X_j^{m+1}X_k^{l+1}}\right|I[X_i<\delta_nnp_n\ \cup\ X_j<\delta_nnp_n\cup X_k<\delta_nnp_n]\right]\\  \nonumber
&\leq&\sum_{i<j<k}\mathbb{E}\left[\left|\Delta_{ijk}(d_i-\mu_i)^{r}(d_j-\mu_j)^{m}(d_k-\mu_k)^{l}\right|I[X_i<\delta_nnp_n]\right]\\  \nonumber
&&+\sum_{i<j<k}\mathbb{E}\left[\left|\Delta_{ijk}(d_i-\mu_i)^{r}(d_j-\mu_j)^{m}(d_k-\mu_k)^{l}\right|I[X_j<\delta_nnp_n]\right]\\ \label{trieq12} 
&&+\sum_{i<j<k}\mathbb{E}\left[\left|\Delta_{ijk}(d_i-\mu_i)^{r}(d_j-\mu_j)^{m}(d_k-\mu_k)^{l}\right|I[X_k<\delta_nnp_n]\right].
\end{eqnarray}
By Lemma \ref{lem1}, we have
\begin{eqnarray} \nonumber
&&\sum_{i<j<k}\mathbb{E}\left[\left|\Delta_{ijk}(d_i-\mu_i)^{r}(d_j-\mu_j)^{m}(d_k-\mu_k)^{l}\right|I[X_i<\delta_nnp_n]\right]\\  \label{trieq13}
&\leq& n^{k_0+3}\mathbb{P}(d_i<\delta_nnp_n)=\exp\left(-\beta np_n(1+o(1))\right).
\end{eqnarray}
Combining (\ref{trieq9})-(\ref{trieq13}), the last two terms of  
(\ref{edeq}) are $o_P\left(\frac{1}{np_n\sqrt{np_n}}+\frac{\sqrt{p_n}}{np_n}\right)$.

Now we consider the third term and the forth term of (\ref{edeq}). Given $t\geq2$ and $r,m,l$ such that $r+m+l=t$, we have

\begin{eqnarray}\nonumber  
&&\sum_{i<j<k}\mathbb{E}\left[\left|\frac{\Delta_{ijk}(d_i-\mu_i)^r(d_j-\mu_j)^m(d_k-\mu_k)^l}{\mu_i^{r+1}\mu_j^{m+1}\mu_k^{l+1}}\right|\right]\\  \nonumber
    &\leq&\sum_{i<j<k}\frac{\mathbb{E}[\Delta_{ijk}]\mathbb{E}[|(d_i-\mu_i)^r|]\mathbb{E}[|(d_j-\mu_j)^m|]\mathbb{E}[|(d_k-\mu_k)^l|]}{\mu_i^{r+1}\mu_j^{m+1}\mu_k^{l+1}}\\
    &=&O\left(\frac{1}{(np_n)^{\frac{t}{2}}}\right).
\end{eqnarray}
If $t\geq4$, then
\begin{eqnarray}\nonumber  
\sum_{i<j<k}\frac{\Delta_{ijk}(d_i-\mu_i)^r(d_j-\mu_j)^m(d_k-\mu_k)^l}{\mu_i^{r+1}\mu_j^{m+1}\mu_k^{l+1}}=o_P\left(\frac{1}{np_n\sqrt{np_n}}+\frac{\sqrt{p_n}}{np_n}\right).
\end{eqnarray}
Next we consider the case $t=2,3$.

Let $t=2$. Without loss of generality, assume $r=2, m=l=0$ or $r=m=1,l=0$.  Consider $r=2, m=l=0$ first. In this case, we have
\begin{eqnarray}\nonumber  
&&\sum_{i<j<k}\frac{\Delta_{ijk}(d_i-\mu_i)^2-\mathbb{E}[\Delta_{ijk}(d_i-\mu_i)^2]}{\mu_i^{3}\mu_j\mu_k}\\  \nonumber
    &=&\sum_{i<j<k,l\neq s}\frac{\Delta_{ijk}\bar{A}_{il}\bar{A}_{is}}{\mu_i^{3}\mu_j\mu_k}+\sum_{i<j<k,l}\frac{\Delta_{ijk}\bar{A}_{il}^2-\mathbb{E}[\Delta_{ijk}\bar{A}_{il}^2]}{\mu_i^{3}\mu_j\mu_k}\\ \label{trieq14}
    &=&\sum_{i<j<k,l\neq s}\frac{\Delta_{ijk}\bar{A}_{il}\bar{A}_{is}}{\mu_i^{3}\mu_j\mu_k}+\sum_{i<j<k,l}\frac{\Delta_{ijk}(\bar{A}_{il}^2-\mathbb{E}[\bar{A}_{il}^2])}{\mu_i^{3}\mu_j\mu_k}+\sum_{i<j<k,l}\frac{(\Delta_{ijk}-\mathbb{E}[\Delta_{ijk}])\mathbb{E}[\bar{A}_{il}^2]}{\mu_i^{3}\mu_j\mu_k}.
\end{eqnarray}
It is easy to verify that
\begin{eqnarray}\nonumber  
    \mathbb{E}\left[\left(\sum_{i<j<k,l\neq s}\frac{\Delta_{ijk}\bar{A}_{il}\bar{A}_{is}}{\mu_i^{3}\mu_j\mu_k}\right)^2\right]&=&O\left(\frac{p_n}{(np_n)^{3}}\right),\\
        \mathbb{E}\left[\left(\sum_{i<j<k,l}\frac{\Delta_{ijk}(\bar{A}_{il}^2-\mathbb{E}[\bar{A}_{il}^2])}{\mu_i^{3}\mu_j\mu_k}\right)^2\right]&=&O\left(\frac{p_n}{(np_n)^{4}}\right),\\
        \mathbb{E}\left[\left(\sum_{i<j<k,l}\frac{(\Delta_{ijk}-\mathbb{E}[\Delta_{ijk}])\mathbb{E}[\bar{A}_{il}^2]}{\mu_i^{3}\mu_j\mu_k}\right)^2\right]&=&O\left(\frac{p_n}{(np_n)^{4}}\right),
\end{eqnarray}
Hence
\begin{eqnarray}\label{trieq15} 
\sum_{i<j<k}\frac{\Delta_{ijk}(d_i-\mu_i)^2-\mathbb{E}[\Delta_{ijk}(d_i-\mu_i)^2]}{\mu_i^{3}\mu_j\mu_k}=o_P\left(\frac{1}{np_n\sqrt{np_n}}+\frac{\sqrt{p_n}}{np_n}\right).
\end{eqnarray}

Consider $t=2$ and $r=m=1,l=0$. In this case,
\begin{eqnarray}\nonumber  
\sum_{i<j<k}\frac{\Delta_{ijk}(d_i-\mu_i)(d_j-\mu_j)}{\mu_i^{2}\mu_j^2\mu_k}=\sum_{i<j<k,l\neq s}\frac{\Delta_{ijk}\bar{A}_{il}\bar{A}_{js}}{\mu_i^{2}\mu_j^2\mu_k}+\sum_{i<j<k,l}\frac{\Delta_{ijk}\bar{A}_{il}\bar{A}_{jl}}{\mu_i^{2}\mu_j^2\mu_k}.
\end{eqnarray}
It is easy to verify that
\begin{eqnarray*}\nonumber  
    \mathbb{E}\left[\left(\sum_{i<j<k,l\neq s}\frac{\Delta_{ijk}\bar{A}_{il}\bar{A}_{js}}{\mu_i^{2}\mu_j^2\mu_k}\right)^2\right]&=&O\left(\frac{p_n^2}{(np_n)^{4}}\right),\\      \mathbb{E}\left[\left(\sum_{i<j<k,l}\frac{\Delta_{ijk}\bar{A}_{il}\bar{A}_{jl}}{\mu_i^{2}\mu_j^2\mu_k}\right)^2\right]&=&O\left(\frac{p_n^2}{(np_n)^{5}}\right),
\end{eqnarray*}
Hence
\begin{eqnarray}\label{trieq16} 
\sum_{i<j<k}\frac{\Delta_{ijk}(d_i-\mu_i)(d_j-\mu_j)}{\mu_i^{2}\mu_j^2\mu_k}=o_P\left(\frac{1}{np_n\sqrt{np_n}}+\frac{\sqrt{p_n}}{np_n}\right).
\end{eqnarray}

Let $t=3$. Without loss of generality, assume $r=3,m=l=0$ or $r=2, m=1, l=0$ or $r=m=l=1$.

Consider $r=m=l=1$. In this case, one has
\begin{eqnarray*}  
&&\sum_{i<j<k}\frac{\Delta_{ijk}(d_i-\mu_i)(d_j-\mu_j)(d_k-\mu_k)}{\mu_i^{2}\mu_j^2\mu_k^2}\\
&=& \sum_{\substack{i<j<k\\ l\neq s\neq m}}\frac{\Delta_{ijk}\bar{A}_{il}\bar{A}_{jm}\bar{A}_{ks}}{\mu_i^{2}\mu_j^2\mu_k^2}+ 3\sum_{\substack{i<j<k\\ l\neq s}}\frac{\Delta_{ijk}\bar{A}_{il}\bar{A}_{jl}\bar{A}_{ks}}{\mu_i^{2}\mu_j^2\mu_k^2}+ \sum_{\substack{i<j<k\\ l}}\frac{\Delta_{ijk}\bar{A}_{il}\bar{A}_{jl}\bar{A}_{kl}}{\mu_i^{2}\mu_j^2\mu_k^2}
\end{eqnarray*}
It is easy to verify that
\begin{eqnarray*} 
\mathbb{E}\left[\left(\sum_{\substack{i<j<k\\ l\neq s\neq m}}\frac{\Delta_{ijk}\bar{A}_{il}\bar{A}_{jm}\bar{A}_{ks}}{\mu_i^{2}\mu_j^2\mu_k^2}\right)^2\right]&=&O\left(\frac{p_n^3}{(np_n)^{6}}\right),\\
\mathbb{E}\left[\left(\sum_{\substack{i<j<k\\ l\neq s}}\frac{\Delta_{ijk}\bar{A}_{il}\bar{A}_{jl}\bar{A}_{ks}}{\mu_i^{2}\mu_j^2\mu_k^2}\right)^2\right]&=&O\left(\frac{p_n}{(np_n)^{7}}\right),\\
\mathbb{E}\left[\left(\sum_{\substack{i<j<k\\ l}}\frac{\Delta_{ijk}\bar{A}_{il}\bar{A}_{jl}\bar{A}_{kl}}{\mu_i^{2}\mu_j^2\mu_k^2}\right)^2\right]&=&O\left(\frac{p_n^2}{(np_n)^{8}}\right).
\end{eqnarray*}

Hence,
\begin{eqnarray}\label{trieq17} 
\sum_{i<j<k}\frac{\Delta_{ijk}(d_i-\mu_i)(d_j-\mu_j)(d_k-\mu_k)}{\mu_i^{2}\mu_j^2\mu_k^2}=o_P\left(\frac{1}{np_n\sqrt{np_n}}+\frac{\sqrt{p_n}}{np_n}\right).
\end{eqnarray}

Consider $r=2, m=1, l=0$. In this case, one has
\begin{eqnarray*}  
\sum_{i<j<k}\frac{\Delta_{ijk}(d_i-\mu_i)^2(d_j-\mu_j)}{\mu_i^{3}\mu_j^2\mu_k}
&=& \sum_{\substack{i<j<k\\ l\neq s\neq m}}\frac{\Delta_{ijk}\bar{A}_{il}\bar{A}_{im}\bar{A}_{js}}{\mu_i^{3}\mu_j^2\mu_k}+ \sum_{\substack{i<j<k\\ l\neq s}}\frac{\Delta_{ijk}\bar{A}_{il}^2\bar{A}_{js}}{\mu_i^{3}\mu_j^2\mu_k}\\
&&+ 2\sum_{\substack{i<j<k\\ l\neq s}}\frac{\Delta_{ijk}\bar{A}_{il}\bar{A}_{is}\bar{A}_{js}}{\mu_i^{3}\mu_j^2\mu_k}+ \sum_{\substack{i<j<k\\ l}}\frac{\Delta_{ijk}\bar{A}_{il}^2\bar{A}_{jl}}{\mu_i^{3}\mu_j^2\mu_k}.
\end{eqnarray*}

It is easy to verify that
\begin{eqnarray*} 
\mathbb{E}\left[\left(\sum_{\substack{i<j<k\\ l\neq s\neq m}}\frac{\Delta_{ijk}\bar{A}_{il}\bar{A}_{im}\bar{A}_{js}}{\mu_i^{3}\mu_j^2\mu_k}\right)^2\right]&=&O\left(\frac{p_n}{(np_n)^{5}}\right),\\
\mathbb{E}\left[\left(\sum_{\substack{i<j<k\\ l\neq s}}\frac{\Delta_{ijk}\bar{A}_{il}^2\bar{A}_{js}}{\mu_i^{3}\mu_j^2\mu_k}\right)^2\right]&=&O\left(\frac{p_n}{(np_n)^{4}}\right),\\
\mathbb{E}\left[\left(\sum_{\substack{i<j<k\\ l\neq s}}\frac{\Delta_{ijk}\bar{A}_{il}\bar{A}_{is}\bar{A}_{js}}{\mu_i^{3}\mu_j^2\mu_k}\right)^2\right]&=&O\left(\frac{p_n^2}{(np_n)^{6}}\right),\\
\mathbb{E}\left[\left|\sum_{\substack{i<j<k\\ l}}\frac{\Delta_{ijk}\bar{A}_{il}^2\bar{A}_{jl}}{\mu_i^{3}\mu_j^2\mu_k}\right|\right]&=&O\left(\frac{p_n}{(np_n)^{2}}\right).
\end{eqnarray*}

Hence
\begin{eqnarray} \label{trieq18}
\sum_{i<j<k}\frac{\Delta_{ijk}(d_i-\mu_i)^2(d_j-\mu_j)}{\mu_i^{3}\mu_j^2\mu_k}
=o_P\left(\frac{1}{np_n\sqrt{np_n}}+\frac{\sqrt{p_n}}{np_n}\right).
\end{eqnarray}

Consider $r=3, m=l=0$. In this case, one has
\begin{eqnarray*}  
&&\sum_{i<j<k}\frac{\Delta_{ijk}(d_i-\mu_i)^3-\mathbb{E}[\Delta_{ijk}(d_i-\mu_i)^3]}{\mu_i^{4}\mu_j\mu_k}\\
&=& \sum_{\substack{i<j<k\\ l\neq s\neq m}}\frac{\Delta_{ijk}\bar{A}_{il}\bar{A}_{im}\bar{A}_{is}}{\mu_i^{4}\mu_j\mu_k}+ 3\sum_{\substack{i<j<k\\ l\neq s}}\frac{\Delta_{ijk}\bar{A}_{il}^2\bar{A}_{is}}{\mu_i^{4}\mu_j\mu_k}\\
&&+ \sum_{\substack{i<j<k\\ l}}\frac{\Delta_{ijk}(\bar{A}_{il}^3-\mathbb{E}[\bar{A}_{il}^3])}{\mu_i^{4}\mu_j\mu_k}+\sum_{\substack{i<j<k\\ l}}\frac{(\Delta_{ijk}-\mathbb{E}[\Delta_{ijk}])\mathbb{E}[\bar{A}_{il}^3]}{\mu_i^{4}\mu_j\mu_k}
\end{eqnarray*}
It is easy to verify that
\begin{eqnarray*}  
\mathbb{E}\left[\left(\sum_{\substack{i<j<k\\ l\neq s\neq m}}\frac{\Delta_{ijk}\bar{A}_{il}\bar{A}_{im}\bar{A}_{is}}{\mu_i^{4}\mu_j\mu_k}\right)^2\right]&=&O\left(\frac{p_n}{(np_n)^4}\right),\\
\mathbb{E}\left[\left(\sum_{\substack{i<j<k\\ l\neq s}}\frac{\Delta_{ijk}\bar{A}_{il}^2\bar{A}_{is}}{\mu_i^{4}\mu_j\mu_k}\right)^2\right]&=&O\left(\frac{p_n}{(np_n)^4}\right),\\
\mathbb{E}\left[\left|\sum_{\substack{i<j<k\\ l}}\frac{\Delta_{ijk}(\bar{A}_{il}^3-\mathbb{E}[\bar{A}_{il}^3])}{\mu_i^{4}\mu_j\mu_k}\right|\right]&=&O\left(\frac{1}{(np_n)^2}\right),\\
\mathbb{E}\left[\left|\sum_{\substack{i<j<k\\ l}}\frac{(\Delta_{ijk}-\mathbb{E}[\Delta_{ijk}])\mathbb{E}[\bar{A}_{il}^3]}{\mu_i^{4}\mu_j\mu_k}\right|\right]&=&O\left(\frac{1}{(np_n)^2}\right).
\end{eqnarray*}
Hence,
\begin{eqnarray*}  
\sum_{i<j<k}\frac{\Delta_{ijk}(d_i-\mu_i)^3-\mathbb{E}[\Delta_{ijk}(d_i-\mu_i)^3]}{\mu_i^{4}\mu_j\mu_k}
=o_P\left(\frac{1}{np_n\sqrt{np_n}}+\frac{\sqrt{p_n}}{np_n}\right).
\end{eqnarray*}

Note that
\begin{eqnarray*}
\sum_{i<j<k}\frac{\bar{A}_{ij}\mu_{jk}\mu_{ki}+\mu_{ij}\bar{A}_{jk}\mu_{ki}+\mu_{ij}\mu_{jk}\bar{A}_{ki}}{\mu_i\mu_j\mu_k}
    &=&\frac{1}{6}\sum_{i\neq j\neq k}\frac{\bar{A}_{ij}\mu_{jk}\mu_{ki}+\mu_{ij}\bar{A}_{jk}\mu_{ki}+\mu_{ij}\mu_{jk}\bar{A}_{ki}}{\mu_i\mu_j\mu_k}
\end{eqnarray*}

Let $\gamma_{ij}=\sum_{k\not\in\{i,j\}}\frac{\mu_{jk}\mu_{ki}}{\mu_i\mu_j\mu_k}$.  Then
\begin{eqnarray*}
\sum_{i\neq j\neq k}\frac{\bar{A}_{ij}\mu_{jk}\mu_{ki}}{\mu_i\mu_j\mu_k}
    =\sum_{i\neq j}\bar{A}_{ij}\sum_{k\not\in\{i,j\}}\frac{\mu_{jk}\mu_{ki}}{\mu_i\mu_j\mu_k}
    =2\sum_{i<j}\gamma_{ij}\bar{A}_{ij}.
\end{eqnarray*}
Then
\begin{eqnarray*}
\sum_{i<j<k}\frac{\bar{A}_{ij}\mu_{jk}\mu_{ki}+\mu_{ij}\bar{A}_{jk}\mu_{ki}+\mu_{ij}\mu_{jk}\bar{A}_{ki}}{\mu_i\mu_j\mu_k}
    &=&\sum_{i<j}\gamma_{ij}\bar{A}_{ij}.
\end{eqnarray*}

\begin{eqnarray*}
    &&\sum_{i< j< k, l\notin\{i,j,k\}}\left(\frac{\mu_{ij}\mu_{jk}\mu_{ki}\bar{A}_{il}}{\mu_i^2\mu_j\mu_k}+\frac{\mu_{ij}\mu_{jk}\mu_{ki}\bar{A}_{jl}}{\mu_i\mu_j^2\mu_k}+\frac{\mu_{ij}\mu_{jk}\mu_{ki}\bar{A}_{kl}}{\mu_i\mu_j\mu_k^2}\right)\\
    &=&\frac{1}{6}\sum_{i\neq j\neq k\neq l}\left(\frac{\mu_{ij}\mu_{jk}\mu_{ki}\bar{A}_{il}}{\mu_i^2\mu_j\mu_k}+\frac{\mu_{ij}\mu_{jk}\mu_{ki}\bar{A}_{jl}}{\mu_i\mu_j^2\mu_k}+\frac{\mu_{ij}\mu_{jk}\mu_{ki}\bar{A}_{kl}}{\mu_i\mu_j\mu_k^2}\right)
\end{eqnarray*}

Let $\eta_{i}=\sum_{\substack{j\neq k}}\frac{\mu_{ij}\mu_{jk}\mu_{ki}}{\mu_i^2\mu_j\mu_k}$. Then
\begin{eqnarray*}
\sum_{i\neq j\neq k\neq l}\frac{\mu_{ij}\mu_{jk}\mu_{ki}\bar{A}_{il}}{\mu_i^2\mu_j\mu_k}
  =\sum_{i\neq l}\bar{A}_{il}\sum_{\substack{j\neq k\\ j,k\not\in\{i,l\}}}\frac{\mu_{ij}\mu_{jk}\mu_{ki}}{\mu_i^2\mu_j\mu_k}=\sum_{i\neq l}\eta_{i}\bar{A}_{il}+\sum_{i\neq l}r_n\bar{A}_{il}, 
\end{eqnarray*}
where $r_n=O\left(\frac{p_n^3}{(np_n)^4}\right)$.
Hence
\begin{eqnarray*}
&&\sum_{i< j< k, l\notin\{i,j,k\}}\left(\frac{\mu_{ij}\mu_{jk}\mu_{ki}\bar{A}_{il}}{\mu_i^2\mu_j\mu_k}+\frac{\mu_{ij}\mu_{jk}\mu_{ki}\bar{A}_{jl}}{\mu_i\mu_j^2\mu_k}+\frac{\mu_{ij}\mu_{jk}\mu_{ki}\bar{A}_{kl}}{\mu_i\mu_j\mu_k^2}\right)\\
    &=&\frac{1}{2}\sum_{i\neq l}\eta_{i}\bar{A}_{il}+\frac{1}{2}\sum_{i\neq l}r_n\bar{A}_{il}\\
    &=&\frac{1}{2}\sum_{i<j}(\eta_{i}+\eta_{j})\bar{A}_{ij}+O_P\left(\frac{p_n^2}{(np_n)^3}\right).
\end{eqnarray*}

Then we get

If $\alpha>\frac{1}{2}$, then 
\begin{eqnarray*} 
    \mathcal{T}_n-\mathbb{E}[\mathcal{T}_n]
=\sum_{i<j<k}\frac{\bar{A}_{ij}\bar{A}_{jk}\bar{A}_{ki}}{\mu_i\mu_j\mu_k}+o_P\left(\frac{1}{np_n\sqrt{np_n}}\right).
\end{eqnarray*}

If $\alpha<\frac{1}{2}$, then 
\begin{eqnarray*}
     \mathcal{T}_n-\mathbb{E}[\mathcal{T}_n]   &=&\sum_{i<j}\left(\gamma_{ij}-\frac{\eta_i+\eta_j}{2}\right)\bar{A}_{ij}
    +o_P\left(\frac{1}{n\sqrt{p_n}}\right).
\end{eqnarray*}

If $\alpha=\frac{1}{2}$, then 
\begin{eqnarray*}\nonumber
    \mathcal{T}_n-\mathbb{E}[\mathcal{T}_n]  &=&\sum_{i<j<k}\frac{\bar{A}_{ij}\bar{A}_{jk}\bar{A}_{ki}}{\mu_i\mu_j\mu_k}+\sum_{i<j}\left(\gamma_{ij}-\frac{\eta_i+\eta_j}{2}\right)\bar{A}_{ij}+o_P\left(\frac{1}{n\sqrt{p_n}}+\frac{1}{np_n\sqrt{np_n}}\right).
\end{eqnarray*}

By Lemma \ref{nlem1}, the proof is complete.

\qed

\subsection{Proof of Corollary \ref{rank1cor}}
We only need to find the orders of $v_{1n},v_{2n}$ and show $v_{2n}=o(v_{1n})$. Note that $w=\Theta(n)$. Given distinct indices $i,j,k$, we have
\begin{eqnarray*}
    \mu_i&=&p_nw_i\sum_{l\neq i }w_l=wp_nw_i\left(1+O\left(\frac{1}{n}\right)\right),\\
    v_{1n}^2&=&\frac{n^3}{6w^6p_n^3}(1+o(1)).
\end{eqnarray*}

Moreover, direct calculation yields
\begin{eqnarray*}
\eta_{i}&=&\sum_{\substack{j\neq k}}\frac{\mu_{ij}\mu_{jk}\mu_{ki}}{\mu_i^2\mu_j\mu_k}\\
&=&\sum_{\substack{j\neq k}}\frac{p_n^3w_i^2w_j^2w_k^2}{w^4p_n^4w_i^2w_jw_k\left(1+O\left(\frac{1}{n}\right)\right)}\\
&=&\frac{1}{w^4p_n}\sum_{\substack{j\neq k}}w_jw_k+O\left(\frac{1}{n^3p_n}\right)\\
&=&\frac{1}{w^2p_n}+O\left(\frac{1}{n^3p_n}\right),
\end{eqnarray*}
and
\begin{eqnarray*}
\gamma_{ij}&=&\sum_{k\not\in\{i,j\}}\frac{\mu_{jk}\mu_{ki}}{\mu_i\mu_j\mu_k}\\
&=&\sum_{k\not\in\{i,j\}}\frac{p_n^2w_iw_jw_k^2}{w^3p_n^3w_iw_jw_k\left(1+O\left(\frac{1}{n}\right)\right)}\\
&=&\frac{1}{w^3p_n}\sum_{k\not\in\{i,j\}}w_k+O\left(\frac{1}{n^3p_n}\right)\\
&=&\frac{1}{w^2p_n}+O\left(\frac{1}{n^3p_n}\right).
\end{eqnarray*}

Hence we get
\begin{eqnarray*}
v_{2n}^2&=&\sum_{i<j}\left(\gamma_{ij}-\frac{\eta_i+\eta_j}{2}\right)^2\mu_{ij}(1-\mu_{ij})\\
&=&O\left(\frac{1}{n^6p_n^2}\right)\sum_{i<j}\mu_{ij}(1-\mu_{ij})\\
&=&O\left(\frac{1}{n^4p_n}\right).
\end{eqnarray*}
Then $v_{2n}=o(v_{1n})$. Then the proof is complete.

\qed


\begin{thebibliography}{9}


\bibitem{ARS01}
Agiwal, M., Roy, A., Saxena, N. (2016). Next generation 5G wireless networks: A comprehensive survey. \textit{IEEE communications surveys and tutorials}, 18(3), 1617-1655.


\bibitem{BW00}
Barrat, A. and M. Weigt (2000). On the properties of small-world network models. \textit{The
European Physical Journal B-Condensed Matter and Complex Systems} 13 (3), 547–560.


\bibitem{BM06}
Bianconi, G. and Marsili, M. (2006).
Number of cliques in random scale-free network ensembles,\textit{Physica D: Nonlinear Phenomena}, 224,:1-6.


\bibitem{BBA01}
Becq, G. J. C., Barbier, E. L., Achard, S. (2020). Brain networks of rats under anesthesia using resting-state fMRI: comparison with dead rats, random noise and generative models of networks. Journal of Neural Engineering, 17(4), 045012.

\bibitem{BCH20}
Bogerd,K., Castro, R., and  Hofstad, R.(2020).
Cliques in rank-1 random graphs: The role of inhomogeneity,\textit{Bernoulli}, 26(1): 253-285 .






\bibitem{CLX01}
Cao, L., Li, L., Huang, Z., Xia, F., Huang, R., Ma, Y., Ren, Z. (2023). Functional network segregation is associated with higher functional connectivity in endurance runners. Neuroscience Letters, 812, 137401

\bibitem{CGL16}
 Chiasserini, C.F.,  Garetto, M. and  Leonardi, E. (2016). Social Network De-Anonymization Under
Scale-Free User Relations. \textit{IEEE/ACM Transactions on Networking} 24 (6):3756–3769.

\bibitem{CV1}
Costanzo, M., VanderSluis, B., Koch, E. N., Baryshnikova, A., Pons, C., Tan, G., Boone, C. (2016). A global genetic interaction network maps a wiring diagram of cellular function. Science, 353(6306), aaf1420.





\bibitem{CHHS21}
Chakrabarty, A., Hazra, S. R., Hollander, F. D. and Sfragara, M.(2021).
Spectra of adjacency and Laplacian matrices of inhomogeneous  Erd\"{o}s-R\'{e}nyi random graphs,
\textit{Random matrices: Theory and applications}, 10(1),215009.

\bibitem{CHHS20}
Chakrabarty, A., Hazra, S. R., Hollander, F. D. and Sfragara, M.(2020).
Large deviation principle for the maximal eigenvalue of
inhomogeneous Erd\"{o}s-R\'{e}nyi random graphs,
\textit{Journal of Theoretical Probability}, https://doi.org/10.1007/s10959-021-01138-w


\bibitem{CCH20}
Chakrabarty, A., Chakrabarty, S. and Hazra, R. S.(2020).
Eigenvalues outside the bulk of  of
inhomogeneous Erd\"{o}s-R\'{e}nyi random graphs,
\textit{Journal of Statistical Physics}, 181: 1746-1780.




\bibitem{HH14} Hall, P. and Heyde, C. C. (2014).
\textit{Martingale limit theory and its application}. Academic press.




\bibitem{MLS22}
Litvak, N., Michielan, R., and Stegehuis, C. (2022).
Detecting hyperbolic geometry in networks: Why triangles are not enough, \textit{Physical Review E} 106 (5), 054303.






\bibitem{LSK1}
Leskovec, J., Singh, A., Kleinberg, J. (2006). Patterns of influence in a recommendation network. In Pacific-Asia Conference on Knowledge Discovery and Data Mining (pp. 380-389). Berlin, Heidelberg: Springer Berlin Heidelberg.

\bibitem{N03}
Newman, M. (2003). The Structure and Function of Complex Networks. \textit{SIAM review} 45, (2), 167–256.


\bibitem{NS21}
Nakajima, K., and K. Shudo.(2021). Measurement error of network clustering coefficients under
randomly missing nodes. \textit{Scientific Reports} 11 (1):1–14.


\bibitem{RFV1}
Radicchi, F., Fortunato, S., Vespignani, A. (2011). Citation networks. Models of science dynamics: Encounters between complexity theory and information sciences, 233-257.

\bibitem{RPW05}
Robins, G., P. Pattison, and J. Woolcock. (2005). Small and other worlds: global network structures from
local processes. \textit{American Journal of Sociology} 110 (4):894–936



\bibitem{TPL01}
Tooley, U. A., Park, A. T., Leonard, J. A., Boroshok, A. L., McDermott, C. L., Tisdall, M. D., Mackey, A. P. (2022). The age of reason: Functional brain network development during childhood. Journal of Neuroscience, 42(44), 8237-8251.

\bibitem{TOSHK06}
Toivonen, R., Onnela, J. Saramaki, J. Hyvonen, J. and Kaski, K.(2006). A model for social networks, \textit{Physica A: Statistical Mechanics and its Applications}, 371(2): 851-860.

\bibitem{VSD01}
Van Wijk, B. C., Stam, C. J., Daffertshofer, A. (2010). Comparing brain networks of different size and connectivity density using graph theory. PloS one, 5(10), e13701.

\bibitem{WS98}
Watts, D. J. and S. H. Strogatz (1998). Collective dynamics of ‘small-world’networks.
\textit{Nature} 393 (6684), 440–442

\bibitem{XKB1}
Xu, T., Cullen, K. R., Mueller, B., Schreiner, M. W., Lim, K. O., Schulz, S. C., Parhi, K. K. (2016). Network analysis of functional brain connectivity in borderline personality disorder using resting-state fMRI. NeuroImage: Clinical, 11, 302-315.



\bibitem{Y23}
Yuan, M.(2023). On the Randić index and its variants of network data. \textit{TEST}, \url{https://doi.org/10.1007/s11749-023-00887-6}



\bibitem{Y23b}
Yuan, M. (2023). Asymptotic distribution of degree-based topological indices, \textit{MATCH Commun. Math. Comput. Chem.}, 91(1): 135-196.


\bibitem{Y23c}
 Yuan, M. (2023). On the Renyi index of random graphs, \textit{Statistical Papers}, Accepted.  \url{https://doi.org/10.1007/s00362-023-01463-8}

\bibitem{ZY23}
Zhao, X. and Yuan, M. (2023). Robustness of clustering coefficients, \textit{Communication in Statistics-Theory and Methods}, Accepted.
 


\end{thebibliography}
\end{document}